\documentclass[11pt]{article}
\usepackage{amssymb}
\usepackage{times}

\title{Classes de Wadge potentielles et th\'eor\`emes d'uniformisation partielle.\indent}
\author{Dominique LECOMTE}
\date{\it ~Fund. Math.~\rm 143 (1993), 231-258}

\newcommand{\Ana}{{\it\Sigma}^{1}_{1}}

\newcommand{\Borel}{{\it\Delta}^{1}_{1}}

\newcommand{\borel}{{\bf\Delta}^{1}_{1}}
\newcommand{\boraone}{{\bf\Sigma}^{0}_{1}}
\newcommand{\boratwo}{{\bf\Sigma}^{0}_{2}}

\newcommand{\boraxi}{{\bf\Sigma}^{0}_{\xi}}
\newcommand{\borxi}{{\bf\Delta}^{0}_{\xi}}
\newcommand{\bortwo}{{\bf\Delta}^{0}_{2}}
\newcommand{\borone}{{\bf\Delta}^{0}_{1}}

\newcommand{\bormone}{{\bf\Pi}^{0}_{1}}
\newcommand{\bormtwo}{{\bf\Pi}^{0}_{2}}

\newcommand{\bormxi}{{\bf\Pi}^{0}_{\xi}}

\newtheorem{thm} {Th\'eor\`eme} [section]
\newtheorem{defi} [thm] {D\'efinition}
\newtheorem{cor} [thm] {Corollaire}
\newtheorem{lem} [thm] {Lemme}
\newtheorem{rem} [thm] {Remarque}
\newtheorem{rems} [thm] {Remarques}
\newtheorem{prop} [thm] {Proposition}
\newtheorem{rap} [thm] {Rappel}
\newtheorem{propdefi} [thm] {Proposition et d\'efinition}
\newtheorem{exe} [thm] {Exemple}

\begin{document}

\maketitle

\noindent {\footnotesize {\bf R\'esum\'e.} On cherche \`a donner une construction aussi simple que possible d'un bor\'elien donn\'e d'un produit de deux espace polonais. D'o\`u l'introduction de la notion de classe de Wadge potentielle. On \'etudie notamment ce que signifie ``ne pas \^etre potentiellement ferm\'e", en montrant des r\'esultats de type Hurewicz. Ceci nous am\`ene naturellement \`a des th\'eor\`emes d'uniformisation partielle, sur des parties ``grosses", au sens du cardinal ou de la cat\'egorie.}\bigskip\smallskip

\noindent\bf {\Large 0 Introduction.}\rm\bigskip

 La notion de classe de Wadge permet de mesurer la
complexit\'e topologique d'un bor\'elien d'un espace
polonais de dimension 0. On peut se demander si la
complexit\'e d'un bor\'elien donn\'e ne diminue pas en
raffinant la topologie polonaise de l'espace ambiant.
Comme on le verra au tout d\'ebut, la r\'eponse est
positive : tout bor\'elien peut \^etre rendu ouvert.\bigskip

 Mais la question analogue se pose avec des topologies
produit, notamment lors de l'\'etude des relations
d'\'equivalence bor\'eliennes ; d'o\`u l'introduction de
la notion de classe de Wadge potentielle (on d\'etaillera
ce lien dans le chapitre 1). Apr\`es avoir effectu\'e quelques 
rappels et \'etabli quelques propri\'et\'es de base, on donne une 
nouvelle d\'emonstration d'un premier r\'esultat 
d'uniformisation partielle, d\'ej\`a prouv\'e par Przymusi\'nski : 
une condition n\'ecessaire et suffisante pour qu'un bor\'elien 
contienne un graphe de fonction bor\'elienne \`a image non 
d\'enombrable. Ensuite, on verra qu'\`a chaque niveau de complexit\'e on 
peut trouver un bor\'elien dont la complexit\'e ne diminue 
pas.\bigskip

 Suite \`a quoi on cherchera \`a savoir si certains r\'esultats 
vrais pour les classes de Wadge peuvent \^etre adapt\'es aux 
classes de Wadge potentielles. En l'occurrence, il s'agit 
d'une part de voir si cette notion correspond \`a une r\'eduction, 
comme dans le cas des classes de Wadge classiques, et on 
verra que non. D'autre part, on cherchera ensuite \`a savoir si on peut obtenir 
des r\'esultats de type Hurewicz, c'est-\`a-dire : ne pas \^etre 
d'une classe donn\'ee, c'est \^etre au moins aussi compliqu\'e 
que des exemples de r\'ef\'erence n'\'etant pas de cette classe.\bigskip

 On obtiendra des r\'esultats partiels pour les potentiellement 
ferm\'es, et pour les petites classes de Wadge (notamment, on 
caract\'erisera les bor\'eliens potentiellement ferm\'es parmi 
les bor\'eliens \`a coupes d\'enombrables, \`a l'aide d'ensembles 
localement \`a projections ouvertes).\bigskip

 Ce qui nous am\`enera \`a de nouveaux r\'esultats 
d'uniformisation, pour des ensembles \`a coupes maigres et des 
$G_{\delta}$ denses, dont le but est d'obtenir d'autres caract\'erisations que celles \'evoqu\'ees ci-dessus.

\vfill\eject

\section{$\!\!\!\!\!\!$ Notations et rappels.}\indent

 On utilisera les notations standard de la th\'eorie
descriptive des ensembles, qui peuvent \^etre trouv\'ees
dans [Mo]. Par exemple, on notera $D_2(\boraone)$ la classe des
diff\'erences de deux ouverts.\bigskip

 $\Gamma$ d\'esignera une famille
d'ensembles, et $\Gamma \lceil X$ les parties de $X$ qui sont dans
$\Gamma$. Par exemple, $\borel \lceil 2^\omega$
d\'esignera l'ensemble des bor\'eliens de $2^\omega$.\bigskip

 Dans un espace polonais r\'ecursivement pr\'esent\'e,
$\it\Sigma$ d\'esignera la topologie engendr\'ee par les
$\Ana$ (c'est la topologie de Gandy-Harrington),
$\it\Delta$ la topologie engendr\'ee par les $\Borel$, ce qu'on pourra noter
aussi $\prec \Borel \succ$ (cette topologie est polonaise : cf [Lo4]).\bigskip

 Si $X$ est un tel espace, on pose $\Omega_X$:=$\{ x \in X
~/~\omega_1^x = \omega_1^{\mbox{CK}}  \}$. Rappelons (cf [Mo]) que $\Omega_X$ est $\Ana$,
et un $\it\Sigma$-ouvert dense (cf [Lo1]). Les traces des $\Ana$ sur 
$\Omega_X$ sont $\borone \lceil (\Omega_X,{\it\Sigma}\lceil \Omega_X)$ ; en
effet, si $A$ est $\Ana$ contenu dans $\Omega_X$ et $f$ est $\Borel$ telle que
$f(x) \in WO \Leftrightarrow x \notin A$, on a 
$$x \notin A \Leftrightarrow
x \notin \Omega_X\ \mbox{ou}\ [x\in \Omega_X\ \mbox{et}\ \exists \xi <
\omega_1^{\mbox{CK}}\ (f(x) \in WO\ \mbox{et}\ \vert f(x)\vert\leq \xi )].$$ 
L'espace $(\Omega_X,{\it\Sigma}\lceil \Omega_X)$ est donc \`a base d\'enombrable
d'ouverts-ferm\'es, donc m\'etrisable s\'eparable ; on sait (cf [Lo1]) que
c'est un espace fortement $\alpha$-favorable, donc c'est un espace polonais,
de dimension 0 par ce qui pr\'ec\`ede.\bigskip
  
 Si $X$ et $Y$ sont des espaces topologiques, $\Pi_X$ (resp.
$\Pi_Y$) d\'esignera la projection de $X\times Y$ sur $X$ (resp. $Y$). Le
symbole $\delta(C)$ d\'esignera le diam\`etre de $C$ pour une distance qui
rende complet l'espace polonais ambiant.\bigskip

  Par ailleurs, je renvoie le lecteur \`a
[Ku] et [Mo] pour ce qui
concerne les notions de base en topologie, en th\'eorie 
descriptive et en th\'eorie effective, et \`a [W] et 
[Lo2] pour les r\'esultats sur les classes de Wadge, qui seront
rappel\'es en cas de besoin. Rappelons tout de m\^eme certains de ces
r\'esultats.\bigskip

 Soit $P_0$ (resp. $P$) un espace polonais de dimension 0,
 et $A_0$ (resp. $A$) un bor\'elien de $P_0$ (resp. $P$). On
 dira que $A$ est dans la $classe\ de\ Wadge\ engendr\acute ee\ par\ A_0$ (not\'ee $\langle A_0 \rangle$) s'il existe une 
fonction continue $f$ de $P$ dans $P_0$ telle que 
$A = f^{-1}(A_0)$.\bigskip

 Si $\Gamma$ est une famille d'ensembles, on notera 
$\check\Gamma := \{ \check A~/~A \in \Gamma \}$. 
Par exemple, $\check \boraone = \bormone$. On dira que 
$\Gamma$ est auto-duale si $\Gamma = \check \Gamma$.\bigskip

 Soit $\Gamma$ une famille de bor\'eliens d'espaces polonais
 de dimension 0, stable par image r\'eciproque continue. 
On montre que $\Gamma$ est une classe de Wadge non 
auto-duale si et seulement si $\Gamma$ admet un universel. 
Ainsi les classes de Baire additives et multiplicatives 
sont des classes de Wadge non auto-duales, par exemple.

\vfill\eject

 On ordonne les classes de Wadge par l'inclusion ; ce qui 
pr\'ec\`ede montre que cet ordre n'est pas total : deux classes
 de Wadge non auto-duales duales l'une de l'autre sont 
incomparables. On a cependant le th\'eor\`eme de Wadge : si 
$\Gamma_1$ et $\Gamma_2$ sont des classes de Wadge, alors
 $\Gamma_1 \subseteq \Gamma_2$ ou 
$\Gamma_2 \subseteq \check\Gamma_1$. Cet ordre est 
\'egalement bien fond\'e. On sait aussi qu'une classe de Wadge
 non auto-duale a pour successeur une classe de Wadge 
auto-duale, et qu'une classe de Wadge auto-duale a pour 
successeurs deux classes de Wadge non auto-duales duales 
l'une de l'autre.\bigskip

 On montre \'egalement que si $\Gamma$ est non auto-duale, 
$2^\omega$ contient un vrai $\Gamma$, c'est-\`a-dire un 
\'el\'ement de $\Gamma\setminus\check\Gamma$.\bigskip

 Si $F$ est ferm\'e
 dans l'espace polonais $P$ de dimension 0, dire que $A$ 
est dans $\Gamma\lceil F$ \'equivaut \`a affirmer l'existence 
de $B$ dans $\Gamma\lceil P$ tel que $A = B \cap F$.\bigskip

 Enfin, si $\xi$ est un ordinal d\'enombrable non nul, on 
note $\borxi\mbox{-PU}(\Gamma)$ la classe suivante : 
$$\{ \bigcup_{n \in \omega} 
A_n \cap U_n ~/~ (U_n)~\borxi\mbox{-partition,}~
A_n \in \Gamma \}.$$ 
On montre que si $\Gamma$ est une 
classe de Wadge, $\Gamma = \borone\mbox{-PU}(\Gamma)$, et que 
si de plus $\Gamma \not= \{ \emptyset \}$ et 
$\check \Gamma \not= \{ \emptyset \}$, il existe un plus 
grand ordinal d\'enombrable $\xi$ tel que 
$\Gamma = \borxi\mbox{-PU}(\Gamma)$, appel\'e $niveau$ de 
$\Gamma$ (on convient que $\{ \emptyset \}$ et
 $\check\{ \emptyset \}$ sont de niveau 0). 

\begin{rap} (a) [Ku] Si $X$ est un espace polonais
et $(B_n)$ une suite de bor\'eliens de $X$, il existe une
topologie polonaise sur $X$, plus fine que la topologie initiale, rendant les
$B_n$ ouverts.\bigskip 

\noindent (b) [Lo3] Si $(X,\sigma)$ et $Y$ sont des
espaces polonais et $B$ un bor\'elien de $X \times Y$ ayant ses
coupes verticales $\boraxi$, il existe une topologie polonaise $\sigma'$
sur $X$, plus fine que $\sigma$, telle que $B$ soit $\boraxi$ 
dans $(X,\sigma')\times Y$.\end{rap}

\begin{lem} Soit $(X,\sigma)$ un espace polonais ; alors il existe 
une topologie polonaise de dimension 0 sur $X$ plus fine que $\sigma$.\end{lem}

\noindent\bf D\'emonstration.\rm\ Soit $(U_n^0)$ une base de la topologie de $\sigma_0 := \sigma$. 
A l'aide du rappel 1.1.(a), on trouve une topologie $\sigma_1$ rendant ferm\'es les 
$U_n^0$. Soit $(U_n^1)$ une base de $\sigma_1$. On construit comme ceci par 
r\'ecurrence une suite croissante $(\sigma_n)$ de topologies polonaises sur $X$ 
telle que si $(U_n^p)_{n\in\omega}$ est une base de $\sigma_p$, $U_n^p$ est 
ferm\'e de $(X,\sigma_{p+1})$. Posons $S_n := \{ U_p^j,~j \leq n,~p\in\omega \}$; alors 
$\sigma_n = \prec S_n \succ$, et $\bigcap_{n\in\omega} S_n = S_0$, donc 
$\sigma' := \prec \bigcup_{n\in\omega} S_n \succ$ r\'epond au probl\`eme: 
$(X,\sigma')$ est hom\'eomorphe \`a la diagonale de $\Pi_{n\in\omega} (X,\sigma_n)$, 
qui est ferm\'ee dans $\Pi_{n\in\omega} (X,\sigma_n)$.$\hfill\square$

\begin{prop} Si $X$ est un espace polonais et $(B_n)$ une suite 
de bor\'eliens de $X$, il existe une topologie polonaise de 
dimension 0 sur $X$, 
plus fine que la topologie initiale, rendant les $B_n$ ouverts.\end{prop}

\noindent\bf D\'emonstration.\rm\ On applique le rappel 1.1.(a) et le lemme 1.2.$\hfill\square$

\begin{defi} Soient $X$ et $Y$ des espaces polonais, et $A$ un 
bor\'elien de $X \times Y$. Si $\Gamma$ est une classe de Wadge de bor\'eliens, on 
dira que $A$ est $potentiellement\ dans\ \Gamma$ (ce qu'on notera 
$A \in \mbox{pot}(\Gamma))$ s'il existe des topologies polonaises de dimension 0, 
$\sigma$ (sur $X$) et $\tau$ (sur $Y$), plus fines que les topologies initiales, 
telles que $A$, consid\'er\'e comme partie de $(X, \sigma) \times (Y, \tau)$, soit 
dans $\Gamma$.\end{defi}

\vfill\eject

 Dans l'\'etude des relations d'\'equivalence bor\'eliennes, par exemple dans 
[HKL], on \'etudie le pr\'e-ordre qui suit. Si $E$ (resp. $E'$) est une 
relation d'\'equivalence bor\'elienne sur l'espace polonais $X$ (resp. $X'$), on pose
$$E \leq E' \Leftrightarrow\ [\exists f:X\rightarrow X'\ \mbox{bor\'elienne telle que}\ 
x E y \Leftrightarrow f(x) E' f(y)].$$

 La derni\`ere relation peut s'\'ecrire : $E = (f\times f)^{-1} (E')$ ; or si $E'$ 
est dans $\Gamma$ (ou m\^eme si $E'$ est $\mbox{pot}(\Gamma)$) et 
$E = (f\times f)^{-1} (E')$, $E$ est $\mbox{pot}(\Gamma)$ : si $(U_n)$ et $(V_n)$ sont 
des bases des topologies associ\'ees \`a $E'$, on peut rendre les $f^{-1} (U_n)$ et 
les $f^{-1} (V_n)$ ouverts, par la proposition 1.3, ce qui fournit une 
topologie $\sigma$; on a alors $E\in \Gamma \lceil (X,\sigma) \times (X,\sigma)$. 
Ceci motive l'introduction de la notion de classe de Wadge potentielle.

\begin{rems} (a) Si $A$ est un bor\'elien \`a coupes verticales 
$\boraxi$ d'un produit de deux espaces polonais, alors $A$ est 
$\mbox{pot}(\boraxi)$.\end{rems}

 En effet, on applique le rappel 1.1.(b) et le lemme 1.2.\bigskip

\noindent\it (b) [Lo4] Si $\Gamma$ est 
une classe de Wadge et $B$ est $\Borel (\alpha)$ dans 
$\omega^\omega\times \omega^\omega$, alors $B$ est $\mbox{pot}(\Gamma)$ si et 
seulement si $B$, consid\'er\'e comme partie de 
$\omega^\omega\times \omega^\omega$ muni de la topologie 
${\it\Delta}_{\alpha}^2$, est dans $\Gamma$.\rm\bigskip

 Rappelons enfin deux r\'esultats sur les ensembles maigres.

\begin{prop}\it [Lo1] Si $X$ est un 
espace polonais parfait non vide et $A$ un sous-ensemble 
maigre de $X \times X$, on trouve une copie $P$ de 
$2^\omega$ dans $X$ telle que si $x$ et $y$ sont distincts 
dans $P$, $(x,y)$ n'est pas dans $A$.\end{prop}

\begin{cor} Si $X$ et $Y$ sont des espaces 
polonais parfaits non vides et $A$ un sous-ensemble maigre 
de $X \times Y$, on trouve une copie $P$ (resp. $Q$) de 
$2^\omega$ dans $X$ (resp. $Y$) telles que  
$(P \times Q) \cap A = \emptyset$.\end{cor}

\noindent\bf D\'emonstration.\rm\ On peut supposer que $X$, $Y = \omega^\omega$ ; 
en effet, si $(U_n)$ est une base de la topologie de $X$, 
$X' := X \setminus [\bigcup_{n\in\omega} (\overline{U_n} 
\setminus U_n)]$ est $G_{\delta}$ dense de $X$, donc 
polonais parfait, et est de dimension 0. Soit $(x_n)$ une 
suite dense de $X'$ ; alors 
$X'' := X' \setminus \{x_n / n \in \omega \}$ est $G_{\delta}$ 
dense de $X'$, donc polonais de dimension 0, et est 
localement non compact, donc hom\'eomorphe \`a $\omega^\omega$. De 
m\^eme avec $Y$.\bigskip

 On applique alors la proposition 1.6. \`a $X = Y = \omega^\omega$ ; 
ceci fournit une injection continue $\psi$ de $2^\omega$ dans 
$\omega^\omega$ dont l'image est le $P$ de la proposition 1.6 ; 
on peut alors poser $P := \psi''N_{(0)}$ et 
$Q := \psi''N_{(1)}$.$\hfill\square$

\section{$\!\!\!\!\!\!$ Ensembles potentiellement ouverts.}\indent

\begin{rem} Soit $A$ un bor\'elien d'un produit 
de deux espaces polonais. Si $A$ a une de ses projections 
d\'enombrable, donc en particulier si $A$ est d\'enombrable, 
$A$ est $\mbox{pot}(\borone)$.\end{rem}

 En effet, si par exemple $\Pi_X''A$ est d\'enombrable, on 
rend les coupes verticales de $A$ ouvertes-ferm\'ees, par 
la proposition 1.3, de sorte que $A \in \mbox{pot}(\boraone) \cap 
\mbox{pot}(\bormone)$, par la remarque 1.5.(a). Ceci fournit des 
topologies $\sigma$ et $\sigma'$ sur $X$, et $\tau$ et
$\tau'$ sur $Y$, de bases respectives ($U_n$), ($U'_n$),
($V_n$), ($V'_n$).

\vfill\eject

 On remarque ensuite que si $(Z,\mu )$ est polonais et
$\mu'$ est polonaise sur $Z$ et plus fine que $\mu$,
l'application identique,
de $(Z,\mu')$ dans $(Z,\mu )$, est
bijective continue ; son inverse est donc bor\'elienne,
et par cons\'equent, les  bor\'eliens de $(Z,\mu )$ et
$(Z,\mu')$ co\"\i ncident. \bigskip

 On applique alors la proposition 1.3 \`a $X$, ($U_n$) et
($U'_n$), qui sont des bor\'eliens de $X$, d'une part, et
\`a $Y$, ($V_n$) et ($V'_n$) d'autre part, pour avoir le
r\'esultat.$\hfill\square$\bigskip

  Ceci montre en particulier que la notion de classe de
Wadge potentielle n'a d'int\'er\^et que si les espaces
polonais ambiants sont non d\'enombrables.

\begin{prop} Soit $A$ un bor\'elien d'un
produit de deux espaces polonais. Alors $A$ est $\mbox{pot}
(\boraone)$ si et seulement si $A$ est r\'eunion d\'enombrable
de rectangles bor\'eliens.\end{prop}

\noindent\bf D\'emonstration.\rm\ Si $A$ est $\mbox{pot}
(\boraone)$, on utilise la remarque de la preuve
pr\'ec\'edente. Inversement, on applique la
proposition 1.3.$\hfill\square$\bigskip
 
 Il r\'esulte de ceci que tous les bor\'eliens d'un
produit ne peuvent \^etre rendus ouverts en conservant
une topologie produit de topologies polonaises ; par
exemple, la diagonale $\Delta(2^\omega)$ de $2^\omega
\times
 2^\omega$ n'est pas r\'eunion d\'enombrable de rectangles,
donc $\Delta(2^\omega)$ n'est pas $\mbox{pot}(\boraone)$.\bigskip

 On donne maintenant une nouvelle preuve du premier 
th\'eor\`eme d'uniformisation partielle annonc\'e ; assez
curieusement, la discussion s'articule autour de la notion
d'ensembles potentiellement ouverts. Ce th\'eor\`eme sera
en outre appliqu\'e au chapitre 3. Notons que ce th\'eor\`eme 
a d\'ej\`a \'et\'e d\'emontr\'e dans [P].

\begin{thm} Soit $A$ un bor\'elien d'un
produit de deux espaces polonais. Alors $A$ contient un
graphe de fonction injective continue d\'efinie sur une
copie de $2^\omega$ si et seulement si $A$ n'est pas
r\'eunion d\'enombrable de rectangles bor\'eliens dont
l'un des c\^ot\'es est un singleton.\end{thm}

\noindent\bf D\'emonstration.\rm\ Si $A$ contient
un graphe comme dans l'\'enonc\'e, raisonnons par l'absurde :
$$A = \bigcup_{n\in\omega} A_n \times  B_n\mbox{,}$$ 
avec $A_n$ ou $B_n$ singleton ; si $Gr(f) \subseteq A$, $Gr(f)$ s'\'ecrit
$\bigcup_{n\in\omega} [Gr(f \lceil A_n) \cap (X \times 
B_n)]$ ; comme $Gr(f)$ est non d\'enombrable, l'un des
ensembles $Gr(f \lceil A_n) \cap (X \times  B_n)$ est non
d\'enombrable, ce qui contredit l'injectivit\'e de $f$.\bigskip

 Montrons la r\'eciproque.\bigskip

\noindent\bf Premier cas.\rm\ $A$ est $\mbox{pot}(\boraone)$.\bigskip

 Par la proposition 2.2, $A = \bigcup_{n\in\omega} A_n
\times  B_n$, avec $A_n$ et $B_n$ bor\'eliens, donc on peut
trouver $n$ tel que $A_n$ et $B_n$ soient non
d\'enombrables, donc bor\'eliennement isomorphes, disons
par $\varphi$. $A_n$ contient une copie de $2^\omega$, qui
contient un $G_\delta$ dense $G$ sur lequel $\varphi$ est
continue ; $G$ \'etant non d\'enombrable contient une
copie de $2^\omega$, d'o\`u le r\'esultat.

\vfill\eject

\noindent\bf Second cas.\rm\ $A$ est non-$\mbox{pot}(\boraone)$.\bigskip

 Posons $E:=\{ x \in X / A(x)$ est
d\'enombrable$\}$. Si $E$ est co-d\'enombrable, $(E
\times Y) \cap A$ est bor\'elien \`a coupes
d\'enombrables, donc est r\'eunion d\'enombrable de graphes
bor\'eliens (par le th\'eor\`eme de Lusin : cf [Mo]). De
plus, par la remarque 2.1, $(\check E \times Y) \cap A$
est $\mbox{pot}(\borone)$, donc $(E \times Y) \cap A$ n'est pas
$\mbox{pot}(\boraone)$, et l'un des graphes n'est pas $\mbox{pot}(\boraone)$ (qui est stable par r\'eunion d\'enombrable par
la proposition 2.2). Par suite, la fonction correspondante
est \`a image non d\'enombrable, par la remarque 2.1. On
va alors trouver un parfait du domaine de la fonction sur
lequel elle est injective, et on conclut comme au
premier cas.\bigskip

 Si $E$ n'est pas co-d\'enombrable, comme il est
co-analytique, $\check E$ contient une copie de
$2^\omega$ ; il suffit donc de voir que si $X = 2^\omega$
et $A$ est \`a coupes verticales non d\'enombrables, $A$
contient un graphe comme dans l'\'enonc\'e. Posons donc 
$$F:=\{ y \in Y ~/~A(y)\mbox{ est maigre dans }2^\omega \}.$$
Si $F$ est co-d\'enombrable, $(2^\omega \times F) \cap A$ est bor\'elien \`a coupes non d\'enombrables
(donc non vides), donc est uniformisable par une fonction
Baire-mesurable d\'efinie sur $2^\omega$ (par
le th\'eor\`eme de von Neumann). Cette fonction $f$ est
continue sur un $G_\delta$ dense $G$ de $2^\omega$, et
$f''G$ est non d\'enombrable : sinon comme
$G=\bigcup_{\beta \in f''G}
 G \cap f^{-1}(\{\beta\})$, l'un des $G \cap
f^{-1}(\{\beta\})$ serait non maigre et contenu dans
$A(\beta)$, ce qui contredirait le fait que $\beta$ est
dans $F$. On conclut alors comme avant.\bigskip

 Si $F$ n'est pas co-d\'enombrable, comme il est
bor\'elien, $\check F$ contient une copie de
$2^\omega$ ; il suffit donc de voir que si $X = Y =
2^\omega$ et $A$ est bor\'elien \`a coupes non maigres, $A$
contient un graphe comme dans l'\'enonc\'e.\bigskip

 Mais par le th\'eor\`eme de Kuratowski-Ulam, $A$ est non
maigre, donc on peut trouver $s$ et $t$ dans $2^{<\omega}$
telles que $(N_s \times  N_t) \setminus A$ soit maigre. On 
trouve alors, par le corollaire 1.7, deux copies $P$ et 
$Q$ de $2^\omega$ telles que 
$P \times Q \subseteq (N_s \times N_t)\cap A$, et si $\varphi$ 
est un hom\'eomorphisme de $P$ sur $Q$, $Gr(\varphi)$ r\'epond \`a la 
question.$\hfill\square$

\section{$\!\!\!\!\!\!$ Classe de Wadge potentielle d'un bor\'elien.}\indent

 On cherche maintenant \`a diminuer au maximum la
complexit\'e d'un bor\'elien donn\'e d'un produit d'espaces
polonais.

\begin{propdefi} Soit $A$ un bor\'elien d'un produit d'espaces polonais. Alors il existe une
unique classe de Wadge de bor\'eliens, la $classe\ de\ Wadge\ potentielle\ de\ A$, not\'ee $\Gamma_A$, telle que\smallskip

\noindent (i) $A \in \mbox{pot}(\Gamma_A)$,\smallskip

\noindent (ii) Si $\Gamma$ est une classe de Wadge strictement 
contenue dans $\Gamma_A$, alors $A\notin \mbox{pot}(\Gamma)$.\end{propdefi}

\noindent\bf D\'emonstration.\rm\ Montrons
l'existence d'une telle classe, en raisonnant par
l'absurde : si $\Gamma_0$ est la classe de Wadge
$\langle$A$\rangle$ engendr\'ee par $A$, $A$ est $\mbox{pot}
(\Gamma_0)$ et on trouve $\Gamma_1 \subset_{\not=}
\Gamma_0$ telle que $A$ soit $\mbox{pot}(\Gamma_1)$. Par
r\'ecurrence, on construit comme ceci $\Gamma_{n+1}
\subset_{\not=} \Gamma_n$ telles que $A$ soit $\mbox{pot}
(\Gamma_{n+1})$. Mais ceci contredit la bonne fondation de
l'ordre de Wadge.

\vfill\eject

 Montrons l'unicit\'e d'une telle classe : si $\Gamma_1$ et
$\Gamma_2$ v\'erifient (i), (ii) et $\Gamma_1 \not=
\Gamma_2$, $\Gamma_1 \not\subseteq \Gamma_2$ (sinon
$\Gamma_1 \subset_{\not=} \Gamma_2$ et $A$ n'est pas $\mbox{pot}
(\Gamma_1)$), donc $\check\Gamma_2 \subseteq \Gamma_1$, et
de m\^eme $\check\Gamma_1 \subseteq \Gamma_2$, donc
$\Gamma_1 = \check\Gamma_2$ est non auto-duale.\bigskip

 L'ensemble $A$ est dans $\mbox{pot}(\Gamma_1) \cap\mbox{pot}
(\Gamma_2)$, donc comme dans la preuve de la remarque
2.1, on trouve des topologies $\sigma$ et $\tau$ telles
que $A \in (\Gamma_1 \cap \Gamma_2) \lceil (X,\sigma)
\times  (Y,\tau)$. Mais la classe de Wadge $\Gamma$
engendr\'ee par $A$, consid\'er\'e comme partie de
$(X,\sigma) \times (Y,\tau)$, v\'erifie 
$$A \in \mbox{pot}(\Gamma)\ \mbox{et}\ \Gamma \subseteq \Gamma_1 \cap
\check\Gamma_1 \subset_{\not=} \Gamma_1\mbox{,}$$ 
une contradiction.$\hfill\square$\bigskip

 Toute classe de Wadge potentielle est donc une classe de
Wadge ; on peut se demander s'il y a une r\'eciproque.
 L'exemple de la diagonale de $2^\omega$ montre que c'est
vrai pour la classe des ferm\'es, et on va voir que
c'est vrai en g\'en\'eral. On peut m\^eme pr\'eciser ce
r\'esultat, en trouvant un $\mbox{pot}(\Gamma)$
``maximal" ; mais pour ce faire on introduit
une classe de fonctions qui est le candidat naturel pour le
probl\`eme de la r\'eduction \'evoqu\'e dans
l'introduction, comme le montre le lemme suivant.\bigskip

 Soit donc $C_0$ la classe des fonctions telles que
l'image r\'eciproque d'un $\mbox{pot}(\boraone)$ soit $\mbox{pot}
(\boraone)$.

\begin{lem} Soient $X$,
$Y$, $X'$, $Y'$, des espaces polonais, $A$ (resp. $B$) un
bor\'elien de $X\times Y$ (resp. $X'\times Y'$) ; si
$f$, de $X\times Y$ dans $X'\times Y'$, est
dans $C_0$ et r\'eduit $A$ \`a $B$, alors $\Gamma_A$ est
contenue dans $\Gamma_B$.\end{lem}

\noindent\bf D\'emonstration.\rm\ $B \in\mbox{pot}
(\Gamma_B)$, ce qui fournit $\sigma$ et $\tau$ telles
que $B \in \Gamma_B \lceil (X',\sigma) \times
(Y',\tau)$. Soient ($U_n$) et ($V_n$) des bases de $\sigma$
et $\tau$. Comme $f$ est dans $C_0$, par la
proposition 2.2 on trouve des bor\'eliens $U_m^{n,p}$ et
$V_m^{n,p}$ tels que $f^{-1}(U_n \times  V_p) = \bigcup_{m
\in \omega} U_m^{n,p} \times  V_m^{n,p}$. Par la
proposition 1.3, on obtient des topologies $\sigma'$ et
$\tau'$ rendant les $U_m^{n,p}$ et les $V_m^{n,p}$
ouverts, de sorte que $f$, de ${(X,\sigma ') \times
(Y,\tau ')}$ dans ${(X',\sigma) \times (Y',\tau)}$ est
continue. Donc ${A
 \in \Gamma_B \lceil (X,\sigma ') \times
(Y,\tau ')}$ et ${A \in \mbox{pot}(\Gamma_B)}$, d'o\`u le
r\'esultat : sinon ${\Gamma_B \subseteq \check\Gamma_A}$, et
${A \in \mbox{pot}(\Gamma_A) \cap \mbox{pot}(\check \Gamma_A)}$ ; donc
comme dans la preuve de la proposition 3.1, et par abus
de langage, ${A \in \mbox{pot}(\Gamma_A \cap \check \Gamma_A)}$ et
${\Gamma_A = \check\Gamma_A}$ ; d'o\`u ${\Gamma_B
\subset_{\not=} \Gamma_A}$, ce qui contredit  ${A \in\mbox{pot}
(\Gamma_B)}$.$\hfill\square$

\begin{thm} Si
$\Gamma$ est une classe de Wadge, il existe
$B$ dans $\Gamma \lceil \omega^\omega \times  \omega^\omega$ tel que\smallskip

\noindent (i) $\Gamma_B = \Gamma$,\smallskip 

\noindent (ii) $A$ est $\mbox{pot}(\Gamma)\Leftrightarrow 
\exists~f\in C_0$ injective telle que $A = f^{-1}(B)$.\end{thm}

\noindent\bf D\'emonstration.\bigskip

\noindent\bf\ Premier cas.\rm\ $\Gamma$ est non auto-duale.\bigskip
 
 Soit $U$ un universel pour $\Gamma \lceil \omega^\omega$,
$U \subseteq \omega^\omega \times  \omega^\omega$.
Posons $(\alpha)_i(n):=\alpha(2n+i)$, o\`u $i = 0$,
$1$,
$$< \gamma,\beta>(n):=\left\{\!\!
\begin{array}{ll}
\gamma(k)\ \ \mbox{si}\ \ n & \!\!\!\! =2k\mbox{,}\cr & \cr
\beta(k)\ \ \mbox{si}\ \ n & \!\!\!\! =2k+1\mbox{,}
\end{array}\right.$$ 
et $B(\alpha,\beta) \Leftrightarrow
U((\alpha)_0,< (\alpha)_1,\beta >)$.

\vfill\eject

 Alors $B$ est aussi universel pour $\Gamma \lceil \omega^\omega$
: $(.)_1$ est continue, donc si $C$ est dans $\Gamma
\lceil \omega^\omega$, $(.)_1^{-1}(C)$ aussi, et il
existe $\alpha$ dans $\omega^\omega$ tel que $(.)_1^{-1}(C) =
U_\alpha$ ; $<\alpha,0^\omega >$ est donc un
code pour $C$.\bigskip

 $B$ est dans $\Gamma$, donc comme \`a la fin de la
preuve du lemme pr\'ec\'edent, $\Gamma_B \subseteq
\Gamma$. Raisonnons par l'absurde pour montrer (i) :
$\Gamma_B \subset_{\not=} \Gamma$, donc $\Gamma_B \subseteq
\check\Gamma$ ; $B$ est $\mbox{pot}(\Gamma_B)$, ce qui
fournit des topologies $\sigma$ et $\tau$ telles que $B \in
\Gamma_B \lceil (\omega^\omega,\sigma) \times  (\omega^\omega,\tau)$, et
on a $B \in \check\Gamma \lceil (\omega^\omega,\sigma) \times
(\omega^\omega,\tau)$.\bigskip

 L'application identique, de $(\omega^\omega,\sigma)$ dans $\omega^\omega$, est
bijective continue, donc d'inverse bor\'elienne ;
son inverse est donc continue sur un $G_\delta$ dense $G$
de $\omega^\omega$ ; sur $G$, $\sigma$ et la topologie usuelle
co\"\i ncident. $G$ \'etant non d\'enombrable contient une
copie $L$ de $2^\omega$, et comme $\Gamma \not= \check
\Gamma$, on peut trouver $D$ dans 
$(\Gamma \setminus \check\Gamma)\lceil L$ ; $D = E \cap L$, 
o\`u $E \in \Gamma \lceil \omega^\omega$. $B$ \'etant universel, soit 
$\alpha$ dans $\omega^\omega$ tel que $B_{\alpha} = E$. Tout 
comme $B_\alpha$, $E$ est dans 
$\check \Gamma \lceil (\omega^\omega,\sigma)$, donc $E \cap G$ est 
dans $\check \Gamma \lceil G$ car sur $G$, les topologies 
sont identiques. Donc $D$ est dans $\check \Gamma \lceil L$, 
une contradiction qui montre que $\Gamma = \Gamma_B$.\bigskip

 Pour (ii), si $A = f^{-1}(B)$, $A$ est $\mbox{pot}(\Gamma)$ \`a 
cause du lemme pr\'ec\'edent. Inversement, si $A$ est $\mbox{pot}(\Gamma)$, 
on trouve $\sigma'$ et $\tau'$ telles que 
$A \in \Gamma \lceil (X,\sigma') \times (Y,\tau')$. On 
trouve alors des ferm\'es $F$ et $H$ de $\omega^\omega$, et des 
hom\'eomorphismes : $\varphi$, de $(X,\sigma')$ sur $F$, et $\psi$, 
de $(Y,\tau')$ sur $H$. 
$(\varphi \times \psi)''A  \in \Gamma \lceil (F \times H)$, donc 
est la trace sur $F \times H$ de $R \in \Gamma \lceil \omega^\omega 
\times \omega^\omega$. $< .,. >$ est un hom\'eomorphisme, donc il 
existe $\alpha$ dans $\omega^\omega$ tel que $< .,. >''R = U_{\alpha}$, 
ce qui s'\'ecrit $R(\gamma,\beta) \Leftrightarrow 
U(\alpha,<\gamma,\beta>) \Leftrightarrow  B(<\alpha,\gamma>,\beta)$. 
La fonction $f := g ~\circ~ (\varphi \times \psi) ~\circ ~\mbox{Id}$ 
r\'epond \`a la question, si on pose 
$$g: \left\{\!\!
\begin{array}{ll} 
F\times H\!\!\!\! & \rightarrow \omega^\omega \times \omega^\omega\mbox{,}\cr & \cr
(\gamma, \beta) & \mapsto (<\alpha,\gamma>,\beta)\mbox{,}
\end{array}\right.$$ 
puisque par la proposition 2.2 les fonctions de la forme $u \times v$, 
avec $u$ et $v$ bor\'eliennes, sont dans $C_0$.\bigskip

\noindent\bf\ Second cas.\rm\ $\Gamma$ est auto-duale.\bigskip

 On sait qu'alors: ou bien il existe une suite strictement 
croissante $(\Gamma_n)$, cofinale dans $\Gamma$, de 
classes de Wadge non auto-duales telle que 
$$\Gamma = \{ \bigcup_{n\in\omega} A_n \cap U_n~/~(U_n)~ 
 \borone\mbox{-partition,}\ A_n \in \Gamma_n \}\mbox{,}$$
ou bien $\Gamma$ est le successeur d'une classe non 
auto-duale $\Gamma'$ telle que 
$$\Gamma = \{ (A \cap N) \cup (B\setminus N)~/~N\in 
\borone,~A\in \Gamma',~B\in \check\Gamma'\}.$$
Dans la premi\`ere \'eventualit\'e, comme $\Gamma_n$ est non 
auto-duale, on trouve $A_n$ dans $\Gamma_n\lceil \omega^\omega \times \omega^\omega$ 
tel que $\Gamma_{A_n} = \Gamma_n$ et si $B_n$ est dans 
$\Gamma_n \lceil \omega^\omega \times \omega^\omega$, il existe $\alpha_n$ 
dans $\omega^\omega$ tel que $B_n(\gamma,\beta) \Leftrightarrow A_n
(<\alpha_n,\gamma>,\beta)$ (ceci par le premier cas).\bigskip

 Soit 
$$\psi_n : \left\{\!\!
\begin{array}{ll} 
\omega^\omega\!\!\!\! & \rightarrow N_{(n)}\mbox{,}\cr & \cr
\alpha & \mapsto n^{\frown}\alpha\mbox{,}
\end{array}\right.$$ 
et $B := \bigcup_{n\in\omega} (\psi_n \times Id)''A_n$.

\vfill\eject

 La fonction $\psi_n$ \'etant un hom\'eomorphisme, 
$$(\psi_n \times Id)''A_n \in \Gamma_n \lceil N_{(n)} \times 
\omega^\omega \subseteq \Gamma \lceil N_{(n)} \times 
\omega^\omega$$ 
et $B \in \borone\mbox{-PU}(\Gamma) = \Gamma$, donc $\Gamma_B \subseteq \Gamma$.\bigskip

 Si l'inclusion est stricte, il existe $n$ tel que 
$\Gamma_B \subset_{\not=} \Gamma_n = \Gamma_{A_n}$ ; or 
$A_n = (\psi_n \times Id)^{-1}(B)$, et 
$\Gamma_{A_n} \subseteq \Gamma_B$ par le lemme pr\'ec\'edent, 
d'o\`u contradiction. Donc $\Gamma_B = \Gamma$.\bigskip

  Si $A$ est $\mbox{pot}(\Gamma)$, on trouve $\sigma'$, $\tau'$, $F$, 
$H$, $\varphi$, $\psi$, $R$ comme au premier cas. On sait 
qu'on peut trouver une $\borone$-partition $(U_n)$ de 
$\omega^\omega\times \omega^\omega$ et $B_n$ dans $\Gamma_n \lceil 
\omega^\omega\times \omega^\omega$ tels que $R = 
\bigcup_{n\in\omega} U_n \cap B_n$. La fonction 
$f := h ~\circ~ (\varphi \times \psi) ~\circ~\mbox{Id}$ convient, si on 
pose
$$h:\left\{\!\!
\begin{array}{ll} 
F \times H\!\!\!\!  & \rightarrow \omega^\omega\times \omega^\omega\mbox{,}\cr & \cr
(\gamma,\beta) & \mapsto (\psi_n(<\alpha_n,\gamma>),\beta)~\ {\rm si}~\ (\gamma,
\beta)\in U_n\mbox{,}
\end{array}\right.$$ 
puisque si $C$ et $D$ sont des bor\'eliens de $\omega^\omega$, on a 
$$h^{-1}(C\times D) = 
\bigcup_{n\in\omega} U_n \cap (\{\alpha \in \omega^\omega~/~
\psi_n(<\alpha_n,\alpha>)\in C \}\times D) \in \mbox{pot}(\boraone).$$
Dans la seconde \'eventualit\'e, on trouve $A_0$ dans $\Gamma'\lceil 
\omega^\omega\times \omega^\omega$, et $A_1$ dans 
$\check \Gamma'\lceil \omega^\omega\times \omega^\omega$ tels que 
$\Gamma_{A_0} = \Gamma'$, $\Gamma_{A_1} = \check \Gamma'$ ; 
et si $C$ (resp. $D$) est dans $\Gamma'$ (resp. $\check 
\Gamma'$) $\lceil \omega^\omega\times \omega^\omega$, on trouve 
$\alpha_0$ (resp. $\alpha_1$) dans $\omega^\omega$ tels que 
$C(\gamma,\beta) \Leftrightarrow  A_0(<\alpha_0,\gamma>,
\beta)$, $D(\gamma,\beta) \Leftrightarrow A_1(<\alpha_1,
\gamma>,\beta)$.\bigskip

 Si $\varphi_0$ est un hom\'eomorphisme de $\omega^\omega$ sur $\omega^\omega 
\setminus N_{(0)}$, et si $B := (\varphi_0\times Id)''A_1 \cup 
(\psi_0\times Id)''A_0$, $B$ est dans 
$\Gamma \lceil \omega^\omega\times \omega^\omega$, donc $\Gamma_B \subseteq \Gamma$.\bigskip

 Comme $\Gamma$ est le successeur de $\Gamma'$, si 
l'inclusion est stricte, on a $\Gamma_B \subseteq \Gamma'$ 
ou $\Gamma_B \subseteq \check\Gamma'$. Soit par exemple 
$\Gamma_B \subseteq \Gamma'$ ; $A_1 \in \mbox{pot}(\check \Gamma')
 = \mbox{pot}(\Gamma_{A_1})$, et comme $\Gamma'$ est non auto-duale, 
$A_1$ n'est pas $\mbox{pot}(\Gamma')$ car 
$\Gamma_{A_1} = \check\Gamma'$. Donc $A_1$ n'est pas 
$\mbox{pot}(\Gamma_B)$ ; or $A_1 = (\varphi_0 \times Id)^{-1}(B)$, 
donc $\Gamma_{A_1} \subseteq \Gamma_B$, une contradiction. 
La derni\`ere partie est analogue \`a celle de la premi\`ere 
\'eventualit\'e.$\hfill\square$\bigskip

 On cherche maintenant \`a adapter les r\'esultats sur les 
classes de Wadge. Si $C$ (resp. $D$) est bor\'elien de $X$ 
(resp. $Y$), on a
$$\langle C \rangle \subseteq \langle D\rangle\ \Leftrightarrow~\mbox{il existe}~f\ \mbox{continue, de}~X
\ \mbox{dans}~Y\mbox{, telle que}~C = f^{-1}(D).$$ 
Une question analogue se pose pour les classes de Wadge 
potentielles : peut-on trouver une classe de fonctions 
${\cal C}$ telle que si $B$ est bor\'elien de $Z \times T$, on ait
$$\Gamma_A \subseteq \Gamma_B~\Leftrightarrow~\mbox{il existe}~f\ \mbox{dans}~{\cal C}\mbox{, de}~
X \times Y\ \mbox{dans}~ Z \times T\mbox{, telle que}~A = f^{-1}(B).$$ 
Comme on va le voir, la r\'eponse est n\'egative ; la classe 
qui semblait le candidat ``raisonnable", $C_0$ \`a cause du 
lemme pr\'ec\'edent, ne fonctionne pas, et \`a un petit niveau 
(avec $A \in \bormone$ et $B\in \check D_2(\boraone)$).

\vfill\eject

 On notera $\leq_P$ le pr\'e-ordre associ\'e \`a $C_0$ :
$$A \leq_P B\ \Leftrightarrow\ \mbox{il existe}\ f\ \mbox{dans}~C_0\ \mbox{telle que}~A = f^{-1}(B).$$
L'in\'egalit\'e $A \leq_P B$ entra\^\i ne donc l'inclusion 
$\Gamma_A \subseteq \Gamma_B$, par le lemme 3.2. On 
montre maintenant un lemme bien plus fort que n\'ecessaire 
pour introduire le contre-exemple \'evoqu\'e ci-dessus, mais 
qui permettra de mieux comprendre ce qu'on cherche \`a faire 
dans les paragraphes suivants.

\begin{defi} Si $X$ est un espace topologique, 
on dira que $G$, $G_{\delta}$ de $X$, est ${presque\mbox{-}ouvert}$ 
si $G$ est contenu dans l'int\'erieur de son adh\'erence.\end{defi}

\begin{lem} Soient $(C_n)~($resp. $(D_n))$ des 
suites de presque-ouverts non vides de $X$ (resp. $Y$), 
$f_n : C_n \rightarrow D_n$ continues et ouvertes, 
$B := \bigcup_{n\in\omega \setminus \{ 0 \}} Gr(f_n)$, et 
$A$ un bor\'elien de $X \times Y$ contenant $B$; si 
$\overline{B} \setminus A$ contient $Gr(f_0)$, alors $A$ 
est non-$\mbox{pot}(\bormone)$.\end{lem}

\noindent\bf D\'emonstration.\rm\ Sinon, soit $F$ (resp. $G$) un 
$G_{\delta}$ dense de $X$ (resp. $Y$) sur lequel les 
topologies (initiales et fournies par le fait que $A$ soit 
$\mbox{pot}(\bormone)$) co\"\i ncident (on montre leur existence comme 
dans la preuve du th\'eor\`eme 3.3) ; on a $A \cap (F \times G) \in 
\bormone \lceil F \times G$.\bigskip

 Montrons que $Gr(f_n) \subseteq \overline{Gr(f_n) \cap 
(F \times G)}$. Soit $U$ (resp. $V$) un ouvert de $X$ 
(resp. $Y$) tels que $(U \times V) \cap  Gr(f_n) \not= \emptyset$. 
Alors $D_n \cap V \cap G$ est un $G_{\delta}$ dense de 
$D_n \cap V$, donc $f_n^{-1}(V \cap G)$ est un $G_{\delta}$ 
dense de $f_n^{-1}(V)$. Donc $F \cap f_n^{-1}(V)$, puis 
$F \cap f_n^{-1}(V \cap G)$, sont des $G_{\delta}$ denses 
de $f_n^{-1}(V)$ ; ce dernier rencontre donc 
$U \cap f_n^{-1}(V)$ en au moins $\{ x \}$ ; on a alors 
$(x,f_n(x)) \in (U \times V) \cap (F \times G) \cap Gr(f_n)  
\not= \emptyset$.\bigskip

 $Gr(f_0)$ est non vide, donc par ce qui pr\'ec\`ede on trouve 
$(x, y)$ dans $(F\times G)\cap Gr(f_0)$, et on a 
$(x,y) \in (F\times G)\cap \overline{B}\setminus A$ ; il suffit 
donc de voir que $(F\times G)\cap \overline{B} \subseteq 
(F\times G)\cap \overline{B\cap(F\times G)}$. On applique alors le fait 
que $Gr(f_n) \subseteq \overline{Gr(f_n) \cap (F \times G)}$ 
pour avoir la contradiction cherch\'ee.$\hfill\square$

\begin{exe} Soit $D_0 := \{(\alpha,\beta)\in 2^\omega \times 
2^\omega~/~\exists ~!~p \in\omega ~~\alpha(p) \not= \beta(p)\}$. 
Alors $D_0$ est non-$\mbox{pot}(\bormone)$.\end{exe}

 En effet, on applique le lemme pr\'ec\'edent \`a $X = Y = C_n = 
D_n = 2^\omega$, $f_0 (\alpha) = \alpha$, 
$$f_n(\alpha)(p)=\alpha(p) \Leftrightarrow p \not= n-1$$  
si $n>0$, et $A=B$.

\begin{thm} Il n'existe pas de classe de fonctions 
$\cal C$ telle que l'inclusion de $\Gamma_A$ dans $\Gamma_B$ 
soit \'equivalente \`a l'existence de $f$ dans $\cal C$ telle 
que $A = f^{-1}(B)$.\end{thm}

\noindent\bf D\'emonstration.\rm\ Raisonnons par l'absurde ; si $B$ est 
$\mbox{pot}(\boraone)$ et $f$ dans $\cal C$, comme $\Gamma_B \subseteq \boraone$, 
$\Gamma_{f^{-1}(B)} \subseteq \boraone$, donc $f^{-1}(B)$ 
est $\mbox{pot}(\boraone)$. Donc si $\cal C$ existe, $\cal C$ est 
une sous-classe de $C_0$.\bigskip

 Comme on l'a vu avant le lemme 3.2, 
$\Gamma_{\Delta(2^\omega)} = \bormone$, et par 3.6, $\check D_0$ 
est non-$\mbox{pot}(\boraone)$, donc $\bormone \subseteq \Gamma_{\check D_0}$ 
et $\Gamma_{\Delta(2^\omega)} \subseteq  \Gamma_{\check D_0}$. 
Il sufflt donc de voir que $\Delta(2^\omega) \not \leq_P \check D_0$ 
pour avoir la contradiction cherch\'ee.

\vfill\eject

 Raisonnons par l'absurde : il existe $f$ dans $C_0$ telle 
que $\Delta(2^\omega) = f^{-1}(\check D_0)$. Alors 
$f''(\Delta(2^\omega))$ est non d\'enombrable, sinon par la 
remarque 2.1 $f''(\Delta(2^\omega))$ serait $\mbox{pot}(\boraone)$, 
et par suite 
$$\Delta(2^\omega) = f^{-1}(f''(\Delta(2^\omega)))$$ 
aussi, ce qui est exclus.\bigskip

 On peut donc trouver une copie $P$ de $2^\omega$ dans 
$\Delta(2^\omega)$ sur laquelle $f$ est injective ; $f''P$ est 
donc un bor\'elien non d\'enombrable, et ses coupes sont 
d\'enombrables : si par exemple une de ses coupes verticales 
$C$ est non d\'enombrable, soit $(f''P)(\alpha_0)$, 
$\{\alpha_0 \}\times C$ est un rectangle bor\'elien, donc est 
$\mbox{pot}(\boraone)$ ; $f^{-1}(\{\alpha_0 \} \times C)$ est alors 
aussi $\mbox{pot}(\boraone)$, non d\'enombrable car 
$\{\alpha_0 \} \times C \subseteq f''(2^\omega \times 2^\omega)$, 
et contenu dans $\Delta(2^\omega)$, ce qui est contradictoire.\bigskip

 $f''P$ n'est donc pas r\'eunion d\'enombrable de rectangles 
bor\'eliens dont l'un des c\^ot\'es soit un singleton, sinon les 
c\^ot\'es seraient d\'enombrables comme les coupes, et $P$ aussi. 
Par le th\'eor\`eme 2.3, il existe un hom\'eomorphisme $\psi$ 
de $2^\omega$ sur un compact $L$, et une injection continue 
$g$ d\'efinie sur $L$, dont le graphe est contenu dans $f''P$.\bigskip

 Alors si $B := D_0 \cap (L \times g''L)$, $B$ est bor\'elien 
\`a coupes d\'enombrables de $L \times g''L$, donc est maigre 
relativement \`a $L \times g''L$. Donc 
$E := (\psi \times (g\circ \psi))^{-1}(B)$ est maigre 
relativement \`a $2^\omega \times 2^\omega$, et par la 
proposition 1.6, il existe $M$ hom\'eomorphe \`a $2^\omega$ 
tel que si $\alpha$ et $\beta$ sont distincts dans $M$, 
alors $(\alpha,\beta)$ n'est pas dans $E$.\bigskip

 Soit $R := \psi''M \times (g\circ \psi)''M$ ; alors 
$R \subseteq \check D_0$, sinon soit $(\alpha,\beta)$ dans 
$R \cap D_0$ ; $\alpha = \psi(\theta)$ et $\beta = g(\psi(\varepsilon))$, 
o\`u $\theta$, $\varepsilon$ sont dans $M$, et $\theta \not= \varepsilon$, sinon 
$(\alpha,\beta) \in Gr(g) \subseteq f''P \subseteq \check D_0$. 
Donc $(\theta, \varepsilon) \notin E$ et $(\alpha,\beta) \notin B$, 
une contradiction.\bigskip

 De plus, $Gr(g\lceil \psi''M) \subseteq R \cap f''P$, donc 
$R \cap f''P$ est non d\'enombrable, et $f^{-1}(R)$ est 
$\mbox{pot}(\boraone)$ comme $R$, est contenu dans 
$\Delta(2^\omega)$, et est non d\'enombrable, ce qui est 
exclus.$\hfill\square$

\section{$\!\!\!\!\!\!$ R\'esultats de type ``Hurewicz".}\indent

 Dans [Lo-SR], il est d\'emontr\'e le r\'esultat suivant :

\begin{thm} Si $\xi$ est un ordinal d\'enombrable non nul, il existe un
compact $P_\xi$ de dimension 0 et un vrai $\boraxi$ de
$P_\xi$, $A_\xi$, tels que si $A$ est un bor\'elien de
l'espace polonais $X$, on ait : $A \notin\bormxi \lceil X$ si et seulement s'il existe  
$f : P_\xi \rightarrow X$ injective continue telle que $A_\xi = f^{-1}(A)$.\end{thm}

 En fait $P_\xi = 2^\omega$, sauf si $\xi = 1$, auquel
cas $P_1$ est constitu\'e d'une suite convergente
infinie et de sa limite. Ceci implique, avec $B =
f''P_\xi$, que $A$ n'est pas $\bormxi$ si et
seulement s'il existe un ferm\'e $B$ de $X$ tel
que $A\cap B$ soit un vrai $\boraxi$ de $B$. L'ensemble
$A_\xi$ est dit ``test d'Hurewicz".\bigskip

 Dans la suite, on cherchera \`a \'etablir un analogue
\`a ces r\'esultats dans le cas o\`u $\xi = 1$. On y
parviendra partiellement. 

\vfill\eject

 Dans cet esprit, voici la 

\begin{defi} Si
$\Gamma$ est une classe, on dira que
$P_1(\Gamma)$ est v\'erifi\'ee si, pour tout $A$ dans $\Gamma$, 
$A$ n'est pas $\mbox{pot}(\bormone)$ si et seulement s'il existe $B\in \mbox{pot}(\bormone)$ et 
$C\in \mbox{pot}(\boraone)$ tels que $B\cap C = B\cap AÊ\notin \mbox{pot}(\bormone)$.\end{defi}

 En apparence, cette propri\'et\'e n'explique pas ce que
signifie ``$A$ est non-$\mbox{pot}(\bormone)$". Mais
elle ram\`ene le probl\`eme au cas o\`u $A$ est $\mbox{pot}
(D_2(\boraone))$, et on va voir que sous l'hypoth\`ese
``$A$ est $\mbox{pot}(F_\sigma)$", on sait
caract\'eriser quand $A$ est non-$\mbox{pot}(\bormone)$.
Mais il nous faut la 
 
\begin{defi} Si $X$
et $Y$ sont des espaces topologiques, une partie $A$
de $X\times Y$ sera dite $localement$ ${\grave a}\ projections\ ouvertes$  
(ou l.p.o.) si pour tout ouvert $U$ de $X \times Y$, les 
projections de $A\cap U$ sont ouvertes.\end{defi}

 Les ensembles l.p.o. se rencontrent par exemple dans la situation suivante : 
$A$ est $\Ana$ dans un produit de deux espaces polonais r\'ecursivement 
pr\'esent\'es. Si on munit ces deux espaces de leur topologie de 
Gandy-Harrington (celle engendr\'ee par les $\Ana$), $A$ devient l.p.o. dans le 
nouveau produit. C'est essentiellement dans cette situation qu'on utilisera 
cette notion.

\begin{lem} Soient $X$ et $Y$
des espaces polonais, $F$ et $G$ des $G_\delta$ denses de
$X$ et $Y$, et $A$ un $G_\delta$ l.p.o. non vide de
$X\times Y$ ; alors $A\cap(F\times G)$ est non vide.\end{lem}

\noindent\bf D\'emonstration.\rm\ Soient ($U_n$) et ($V_n$) des suites
d'ouverts denses, de $X$ et $Y$, telles que $F =
\bigcap_{n\in\omega} U_n$, $G = \bigcap_{n\in\omega} V_n$,
et $(F_n)$ une suite de ferm\'es de $X\times Y$ telle que
$A = \bigcap_{n\in\omega} \check F_n$. On construit par
r\'ecurrence sur $n$ des suites d'ouverts non vides
($O_n$) et ($T_n$) de $X$ et $Y$ v\'erifiant 
$$\begin{array}{ll} 
& (i)~~~ \delta(O_n),\  \delta(T_n) <2^{-n} \cr 
& (ii)~~ O_n\times T_n \subseteq (U_n\times V_n)\setminus F_n \cr 
& (iii)~ A \cap (O_n\times
T_n) \not= \emptyset \cr 
& (iv)~~ \overline{O_{n+1}}\subseteq O_n,\  \overline{T_{n+1}} \subseteq T_n
\end{array}$$
Admettons avoir construit ces objets ;
($\overline{O_n}$) et ($\overline{T_n}$) sont des suites
d\'ecroissantes de ferm\'es non vides dont les
diam\`etres tendent vers 0, donc leurs intersections sont $\{ x \}$ et $\{ y \}$ ; mais 
$$(x,y) \in O_n \times T_n \subseteq \check F_n \cap (U_n \times V_n)\mbox{,}$$ 
donc $(x,y) \in A \cap (F \times G) \not= \emptyset$.\bigskip

 Admettons avoir construit $(O_p)_{p<n}$ et $(T_p)_{p<n}$ 
v\'erifiant les conditions demand\'ees ; alors par (iii), 
$\Pi_X''(A\cap (O_{n-1}\times T_{n-1}))$ est un ouvert non 
vide de $X$, donc rencontre $U_n$ : l'ensemble 
$$A \cap [(O_{n-1} \cap U_n) \times T_{n-1}]$$ 
est non vide, donc sa projection sur $Y$ est un ouvert non vide de $Y$, 
donc rencontre $V_n$. L'ensemble 
$A \cap [(O_{n-1} \cap U_n) \times (T_{n-1} \cap V_n)]$ est 
non vide, donc contient $(x_n,y_n)$ ; il reste \`a choisir 
deux ouverts $O_n$ et $T_n$ de diam\`etre au plus $2^{-n}$ 
v\'erifiant la double inclusion suivante : 
$$(x_n,y_n) \in O_n \times T_n \subseteq \overline{O_n \times T_n} 
\subseteq [(O_{n-1} \cap U_n) \times (T_{n-1} \cap V_n)] 
\setminus F_n.$$
Ceci termine la d\'emonstration.$\hfill\square$

\vfill\eject

 On remarquera que ce r\'esultat est faux si on ne suppose 
pas que $A$ est un $G_{\delta}$. D\'esignons par $P_{\infty}$ 
l'ensemble des suites de 0 et de 1 comportant une infinit\'e 
de termes \'egaux \`a 1. Si maintenant $A$ d\'esigne 
$(2^\omegaÊ\times 2^\omega) \setminus 
(P_{\infty} \times P_{\infty})$, $A$ est un $K_{\sigma}$ 
qui v\'erifie les autres conditions du lemme et ne rencontre 
pas $P_{\infty} \times P_{\infty}$ !

\begin{thm} (a) Soit $A$ un $\mbox{pot}(F_{\sigma})$ 
d'un produit d'espaces polonais $X \times Y$ ; alors $A$ 
est non-$\mbox{pot}(\bormone)$ si et seulement s'il existe des 
espaces polonais $Z$ et $T$ de dimension 0, une suite $(F_n)$ de 
ferm\'es de $Z \times T$, $f:Z\rightarrow X$ et $g:T\rightarrow Y$ 
injectives continues tels que si $C := \bigcup_{n\in\omega} F_n$, on ait 
$$\begin{array}{ll} 
& (i)~~~\overline{C}\setminus C \not= \emptyset \cr 
& (ii)~~C = (f \times g)^{-1}(A) \cr 
& (iii)~\overline{C} \setminus C~\mbox{et}~F_n~\mbox{sont~l.p.o.}
\end{array}$$
\noindent (b) Soit $A$ un $\mbox{pot}(\boratwo)\cap \mbox{pot}(\bormtwo)$ d'un produit d'espaces 
polonais $X \times Y$ ; alors $A$ est non-$\mbox{pot}(\bormone)$ 
si et seulement s'il existe des espaces polonais $Z$ et $T$ de 
dimension 0, ${C\in \bortwo\lceil Z \times T}$, $f:Z\rightarrow X$ et 
$g:T\rightarrow Y$ injectives continues tels que 
$$\begin{array}{ll} 
& (i)~~~\overline{C} \setminus C \not=\emptyset \cr 
& (ii)~~C = (f\times g)^{-1}(A) \cr 
& (iii)~\overline{C} \setminus C~\mbox{et}~C~\mbox{sont~l.p.o.}
\end{array}$$\end{thm}

\noindent\bf D\'emonstration.\rm\ (a) Raisonnons par
l'absurde : si $A$ est $\mbox{pot}(\bormone)$, $C$ aussi et on
trouve des $G_\delta$ denses $F$ et $G$ de $Z$ et $T$ tels
que $C \cap(F\times G) \in \bormone \lceil F\times G$.\bigskip

 On applique le lemme
pr\'ec\'edent \`a $Z$, $T$, $F$, $G$, et $\overline{C}
\setminus C$, et on a : $[(F\times G) \cap
\overline{C} ]
\setminus C \not= \emptyset$. Il suffit alors de voir que
$[(F\times G)  \cap \overline{C}] \setminus C =
[(F\times G)  \cap \overline{C \cap (F\times G)}]
\setminus C$ pour avoir la contradiction cherch\'ee ; et
il suffit de voir que $F_n \subseteq \overline{F_n \cap
(F\times G)}$.\bigskip

 Si $U$ et $V$ sont des ouverts de $Z$ et $T$ tels que 
$F_n \cap (U\times V)$ soit non vide, on applique le
lemme pr\'ec\'edent \`a $U$, $V$, $F\cap U, G\cap V$,
 et $F_n \cap(U\times V)$ pour voir que $F_n
\cap[(F\cap U)\times (G\cap V)]$ est lui aussi non
vide.\bigskip

 Inversement, on peut supposer que $X$ et $Y$ sont des
ferm\'es de $\omega^\omega$, et pour simplifier l'\'ecriture
qu'ils sont $\Borel$, ainsi que la suite ($G_n$) de
ferm\'es pour ${\it\Delta}^2$ dont $A$ est la r\'eunion (on
applique la remarque 1.5.(b)).\bigskip

 Comme ${\it\Delta}\subseteq {\it\Sigma}$, par une double
application du th\'eor\`eme de s\'eparation on voit que 
$\overline {A}^{{\it\Sigma}^2} = \overline {A}^{{\it\Delta}^2}$ ;
par suite, puisque $\it\Delta$ est polonaise, $\overline
{A}^{{\it\Sigma}^2} \setminus A$ est  un $\Ana$ non vide,
ainsi que $(\overline {A}^{{\it\Sigma}^2} \setminus A)\cap
\Omega_{\omega^\omega\times \omega^\omega}$ ; et puisque
 $\Omega_{\omega^\omega\times\omega^\omega} \subseteq
\Omega_{\omega^\omega}^2$, $(\overline
{A}^{{\it\Sigma}^2} \setminus A)\cap \Omega_{\omega^\omega}^2$ est
lui aussi non vide ; posons 
$$Z := (X\cap \Omega_{\omega^\omega},
{\it\Sigma}\lceil X\cap \Omega_{\omega^\omega})\mbox{,}$$ 
$T := (Y\cap \Omega_{\omega^\omega},{\it\Sigma}\lceil Y\cap \Omega_{\omega^\omega})$, 
$F_n := G_n \cap (Z\times T)$, et prenons pour $f$ et $g$ les applications
identiques.\bigskip

 Alors $(\overline
{A}^{{\it\Sigma}^2} \setminus A)\cap \Omega_{\omega^\omega}^2 = 
\overline{A\cap(Z\times T)}^{{\it\Sigma}^2}
\cap  (Z\times T) \setminus 
(A\cap(Z\times T))$, donc ces objets conviennent. En effet, $\overline{C} \setminus C$ et $F_n$ sont 
$\Ana$ et un ouvert de $Z\times T$ est r\'eunion de rectangles $\Ana$
(les projections des traces de ces rectangles seront donc
$\Ana$).

\vfill\eject

\noindent (b) Si $A$ est non-$\mbox{pot}(\bormone)$, on raisonne comme
dans (a), \`a ceci pr\`es
qu'on pose ${C := A \cap (Z\times T)}$.\bigskip

 La r\'eciproque est analogue \`a celle de (a),
\`a ceci pr\`es que pour montrer l'\'egalit\'e entre 
$$[(F\times G)  \cap \overline{C}] \setminus C$$ 
et $[(F\times G)  \cap \overline{C \cap (F\times G)}]
\setminus C$, il suffit de voir que $C \subseteq 
\overline{C \cap (F\times G)}$ ; si $U$ et $V$ sont des
ouverts de $Z$ et $T$ tels que $C \cap (U\times V)$ soit 
non vide, on applique le lemme 4.4 \`a $U$, $V$, $F\cap U$, $G\cap V$, et 
$C\cap (U\times V)$.$\hfill\square$\bigskip

  On introduit maintenant une propri\'et\'e, qui est du type 
Hurewicz au sens de l'introduction ; \`a ceci 
pr\`es que pour comparer la complexit\'e, on n'a pas de 
r\'eduction sur tout l'espace de d\'epart, mais seulement sur 
un ferm\'e.

\begin{defi} Si $\Gamma$ est une classe, on 
dira que $P_2(\Gamma)$ est v\'erifi\'ee si pour tout $A$ dans 
$\Gamma \lceil X\times Y$, $A$ est non-$\mbox{pot}(\bormone)$ si 
et seulement s'il existe des espaces polonais $Z$ et $T$ de dimension 
0, $D$ dans $D_2(\boraone)\lceil Z \times T$, $f:Z\rightarrow X$ et 
$g:T\rightarrow Y$ injectives continues tels que 
$$\begin{array}{ll} 
& (i)~~~\overline{D}\setminus D \not= \emptyset \cr 
& (ii)~~\overline{D} \cap (f \times g)^{-1}(A) = D \cr 
& (iii)~D~\mbox{et}~\overline{D} \setminus D~\mbox{sont~l.p.o.}
\end{array}$$\end{defi}

\begin{prop} $P_1(\Gamma)$ \'equivaut \`a $P_2(\Gamma)$.\end{prop}

\noindent\bf D\'emonstration.\rm\ Supposons $P_1(\Gamma)$ ; si $A\in \Gamma \setminus 
 \mbox{pot}(\bormone)$, soient $B$
et $C$ fournis par $P_1(\Gamma)$, $\sigma$ et $\tau$
rendant $B$ ferm\'e. Alors comme dans la preuve du th\'eor\`eme 4.5.(b)
on trouve $Z$, $T$, et $D$ comme indiqu\'e et des injections continues 
$F$, de $Z$ dans $(X,\sigma)$, et $G$, de $T$ dans $(Y,\tau)$, 
tels que l'on ait l'\'egalit\'e $D = (F\times G)^{-1}(A) \cap (F\times G)^{-1}(B)$ 
(ceci parce que $B\cap C$ est $\mbox{pot}(D_2(\boraone))$). D'o\`u 
$\overline{D} \cap (F\times G)^{-1}(A) = D$, et il
ne reste qu'\`a revenir aux topologies initiales pour
obtenir $f$ et $g$.\bigskip

 Inversement, si $A \in \mbox{pot}(\bormone)$, $D$ aussi, ce
qui est exclus, comme dans la preuve du th\'eor\`eme 4.5.(b).\bigskip

 Supposons maintenant $P_2(\Gamma)$ ; $\mbox{pot}(\bormone)$ \'etant stable par
intersection finie, $A$ est non-$\mbox{pot}(\bormone)$ si $B$ et $C$ existent.\bigskip

 Inversement, si $A$ est non-$\mbox{pot}(\bormone)$, soient $Z$, $T$, $D$, $f$, et $g$ 
 fournis par $P_2(\Gamma)$ ; $D = U\cap \overline {D}$, o\`u $U$ est ouvert de 
 $Z\times T$, donc $B := (f\times g)''\overline {D}$ et $C :=
(f\times g)''U$ r\'epondent au probl\`eme, par
injectivit\'e de $f$ et $g$ : on a $B \cap A = (f\times g)''D = C \cap B$.$\hfill\square$\bigskip

 Cette proposition permet de comprendre pourquoi,
indirectement, la propri\'et\'e $P_1$ permet de mieux
conna\^\i tre les bor\'eliens non potentiellement ferm\'es. On
\'etablit maintenant cette propri\'et\'e pour
certaines familles de bor\'eliens.

\begin{prop} La propri\'et\'e $P_1$ est v\'erifi\'ee par chacune des classes suivantes :\smallskip
 
\noindent (i) Les ensembles potentiellement $\boratwo$ et potentiellement $\bormtwo$.\smallskip

\noindent (ii) Les bor\'eliens \`a coupes~verticales~co-d\'enombrables.\smallskip 

\noindent (iii) Les relations~d'\'equivalence~bor\'eliennes.\end{prop}

\vfill\eject

Avant de d\'emontrer cette proposition, on donne la
d\'efinition suivante :

\begin{defi} Si $\Gamma$ est une classe, 
on dira que $P_3(\Gamma)$ est v\'erifi\'ee si pour tout $A$ 
dans $\Gamma$, $A$ n'est pas $\mbox{pot}(\boraone)$ si et seulement s'il existe 
$B\in \mbox{pot}(\bormone)$ tel que $A \cap B \in \mbox{pot}(\bormone)$, et tel que pour tout $C\in \mbox{pot}(\boraone)$ on 
ait $A \cap B \not= C \cap B$.\end{defi}

 Il est clair, en raison des formules $B\setminus A = 
B\setminus (A \cap B)$ et $B \cap A = B\setminus 
(B\setminus A)$ que $P_1(\Gamma)$ \'equivaut \`a $P_3(\Gamma)$ 
si $\Gamma$ est auto-duale ; mais cette derni\`ere est plus 
maniable quand il est question de r\'eunions d\'enombrables.\bigskip

\noindent\bf D\'emonstration.\rm\ On a d\'emontr\'e le ``si" dans le cas 
g\'en\'eral ; on montre donc la r\'eciproque dans chacun des 
trois cas.\bigskip

\noindent (i) On montre la chose plus pr\'ecise suivante : si $\Gamma$ 
est un contre-exemple minimal (pour l'ordre de Wadge) \`a 
$P_3(\mbox{pot}(\Gamma))$, $\Gamma$ est non auto-duale et est de 
niveau au moins deux.\bigskip

 Si $\Gamma$ est auto-duale, traitons le premier cas de 
l'alternative \'evoqu\'ee dans la preuve du th\'eor\`eme 3.3 
(l'autre cas \'etant plus simple) : on trouve une partition 
$(U_n)$ de $X \times Y$ en $\mbox{pot}(\borone)$ et $A_n$ dans 
$\mbox{pot}(\Gamma_n)$, o\`u $\Gamma_n \subset_{\not=} \Gamma$, avec 
$A \cap U_n = A_n \cap U_n$ ; $A = \bigcup_{n\in\omega} A_n 
\cap U_n$, donc il existe $n$ tel que $A_n \cap U_n$ ne 
soit pas $\mbox{pot}(\boraone)$. Puisque $A_n \cap U_n$ est 
$\mbox{pot}(\Gamma_n)$, par minimalit\'e de $\Gamma$, il existe 
donc $B'$ dans $\mbox{pot}(\bormone)$ tel que $A_n \cap U_n \cap B'$ 
soit $\mbox{pot}(\bormone)$ et pour tout $C$ dans $\mbox{pot}(\boraone)$, 
$A_n \cap U_n \cap B' \not= B' \cap C$. Il reste \`a poser 
$B := B' \cap U_n$.\bigskip

 Si $\Gamma$ est non auto-duale et de niveau au plus un, 
$\Gamma$ est de niveau un car $\Gamma$ contient les ferm\'es, 
donc est de la forme $SD_{\eta}(\Delta,\Gamma_*)$ :\bigskip 
\noindent{$A = E \cup F$}, o\`u 
$E = \bigcup_{\xi<\eta,n\in\omega} A_{\xi,n} \cap V_{\xi,n} \cap
~^c(\bigcup_{\theta<\xi,p\in\omega} V_{\theta,p})$ 
et $F = A_* \cap ~^c(\bigcup_{\xi<\eta,n\in\omega} V_{\xi,n})$, 
$A_{\xi,n}$ \'etant $\mbox{pot}(\Delta)$, $(V_{\xi,n})_{n\in\omega}$ 
\'etant une suite de $\mbox{pot}(\boraone)$ deux \`a deux disjoints, 
et $A_*$ \'etant $\mbox{pot}(\Gamma_*)$, o\`u 
$\check \Gamma_* \not= \Gamma_* \subseteq \Delta$.\bigskip

\noindent\bf Premier cas.\rm\ $E$ est non-$\mbox{pot}(\boraone)$.\bigskip

 Soient $(V_{\xi,n,p})_{p\in\omega}$ une partition de 
$V_{\xi,n}$ en ensembles $\mbox{pot}(\borone)$,
$$\begin{array}{ll}
V_{n,p}^{\xi} & \!\!\! := V_{\xi,n,p} \cap [\bigcup_{\theta<\xi,q\in\omega} 
A_{\theta,q} \cap V_{\theta,q} \cap ~^c(\bigcup_{\rho<\theta,r\in
\omega} V_{\rho,r})]\mbox{, et}\cr & \cr
D_{n,p}^{\xi} & \!\!\! := [A_{\xi,n} \cap V_{\xi,n,p} \cap 
~^c(\bigcup_{\theta<\xi,q\in\omega} V_{\theta,q})] \cup V_{n,p}^{\xi}.
\end{array}$$

 Alors $E = \bigcup_{\xi<\eta;n,p\in\omega} D_{n,p}^{\xi}$, 
donc on trouve $\xi$ minimal tel qu'il existe $n$ et $p$ 
tels que $D_{n,p}^{\xi}$ ne soit pas $\mbox{pot}(\boraone)$. Comme 
$V_{n,p}^{\xi} = V_{\xi,n,p} \cap \bigcup_{\theta<\xi;q,r\in\omega} D_{q,r}^{\theta}$, 
$V_{n,p}^{\xi}$ est $\mbox{pot}(\boraone)$ par la minimalit\'e de $\xi$ ; or 
$$\Delta = \left\{\!\!\!\!\!\!
\begin{array}{ll} 
& \Gamma_0 \cup \check \Gamma_0\mbox{,}\cr 
& \bigcup_{n\in\omega} \Gamma_n\mbox{,}
\end{array}\right.$$ 
o\`u $(\Gamma_n)$ 
est une suite strictement croissante de classes de Wadge de niveau $\not= 1$, donc sauf si\bigskip
 
(1) ``$\Gamma_0 = \{ \emptyset \}$ et $\Delta = \Gamma_0 \cup 
\check \Gamma_0$",\bigskip

\noindent il existe $\Gamma' \subset_{\not=} \Gamma$ 
telle que $D_{n,p}^{\xi}$ soit $\mbox{pot}(\Gamma')$.

\vfill\eject

 Mais si on a (1), $D_{n,p}^{\xi} = G \cup V_{n,p}^{\xi}$, 
o\`u $G := D_{n,p}^{\xi} \setminus V_{n,p}^{\xi}$ est 
$\mbox{pot}(\bormone)$ ; alors si on pose $B' := \check V_{n,p}^{\xi}$, 
$B'$ est $\mbox{pot}(\bormone)$, $D_{n,p}^{\xi}\cap B' = G$ aussi, 
et si $C$ est $\mbox{pot}(\boraone)$ et $C\cap B' = B'\cap D_{n,p}^{\xi}$, 
$$D_{n,p}^{\xi} = (C \cap B') \cup V_{n,p}^{\xi} = C \cup 
V_{n,p}^{\xi}$$ 
est $\mbox{pot}(\boraone)$ ; par cons\'equent $B'$ 
v\'erifie les conclusions de la propri\'et\'e $P_3$ avec $A := 
D_{n,p}^{\xi}$.\bigskip

 Dans l'autre cas, par minimalit\'e de $\Gamma$, on trouve 
aussi un $B'$ v\'erifiant ces conclusions. Il reste \`a poser 
$B := B' \cap V_{\xi,n,p}$, puisque $A \cap B = B' \cap 
D_{n,p}^{\xi}$.\bigskip

\noindent\bf Second cas.\rm\ $E$ est $\mbox{pot}(\boraone)$.\bigskip

 Dans ce cas $F$ est non-$\mbox{pot}(\boraone)$, et sauf si on a 
(1), on trouve $\Gamma' \subset_{\not=} \Gamma$ de niveau 
au moins deux telle que $F$ soit $\mbox{pot}(\Gamma')$, donc telle 
que $A$ soit $\mbox{pot}(\Gamma')$, et l'hypoth\`ese de minimalit\'e 
s'applique. Si on a (1), $F$ est $\mbox{pot}(\bormone)$, ce qu'on 
a trait\'e au premier cas.\bigskip

 Si maintenant $A$ est $\mbox{pot}(\boratwo)\cap \mbox{pot}(\bormtwo)$, il existe un ordinal 
d\'enombrable $\eta >1$ tel que $\check A$ soit $\mbox{pot}(D_{\eta}(\boraone))$, 
par le th\'eor\`eme de Hausdorff. Comme $\Gamma$ est de niveau 
au moins deux, 
$$D_{\eta}(\boraone) \subset_{\not=} \bortwo \subseteq \Gamma .$$
Par ce qui pr\'ec\`ede, il existe $B$ dans $\mbox{pot}(\bormone)$ tel 
que $B\setminus A$ soit $\mbox{pot}(\bormone)$ et pour tout $H$ 
dans $\mbox{pot}(\boraone)$, $B\setminus A \not= B \cap H$. 
Alors $C := A \cup \check B$ r\'epond \`a la question.\bigskip

\noindent (ii) Soit $A$ un bor\'elien \`a coupes verticales 
co-d\'enombrables ; on montre plus que la propri\'et\'e $P_1$ :\bigskip

 (2) \it Il existe des copies $P$ et $Q$ de $2^\omega$, et un 
hom\'eomorphisme $\phi$ de $P$ dans $Q$ tels que 
$$\mbox{Gr}(\phi) = (P \times Q) \setminus A.$$\rm
 Comme dans la d\'emonstration du th\'eor\`eme 2.3, $\check A = 
\bigcup_{n\in\omega} Gr(f_n)$, o\`u les fonctions $f_n$ sont 
bor\'elien-nes et d\'efinies sur des bor\'eliens $B_n$ et on 
trouve un hom\'eomorphisme $\phi_0$ de $2^\omega$ sur $R \subseteq 
B_0$, avec $f_0\lceil R$ injective continue. Posons 
$$E(\alpha,\beta)\ \Leftrightarrow\ \mbox{il existe}\ n > 0\ \mbox{tel que}\ 
\phi_0(\alpha)\in B_n\ \mbox{et}\ f_0(\phi_0(\beta)) = 
f_n(\phi_0(\alpha)).$$
Alors $E$ est bor\'elien \`a coupes verticales d\'enombrables 
de $2^\omega \times 2^\omega$, donc est maigre relativement 
\`a $2^\omega \times 2^\omega$ et par la proposition 1.6 il 
existe dans $2^\omega$ une copie $L$ de $2^\omega$ telle que si 
$\alpha$ et $\beta$ sont distincts dans $L$, $(\alpha,\beta)$ 
n'est pas dans $E$.\bigskip

 Si $\alpha$ et $\beta$ sont distincts dans $\phi_0''L$, 
alors pour tout $n > 0$, $\beta \notin B_n$ ou $f_0(\alpha) 
\not= f_n(\beta)$. Il reste \`a poser $P := \phi_0''L$, 
$Q := (f_0 \circ \phi_0)''L$, $\phi := f_0\lceil P$.\bigskip

 L'ensemble $(P \times Q)\setminus Gr(\phi)$ est non-$\mbox{pot}(\bormone)$, 
 sinon $Gr(\phi)$ serait $\mbox{pot}(\boraone)$ ainsi que 
$$\Delta(2^\omega) = (\mbox{Id}\times \phi)^{-1}(Gr(\phi)).$$
Les ensembles $B := P \times Q$ et $C := \check Gr(\phi)$ 
r\'epondent donc au probl\`eme.

\vfill\eject

\noindent (iii) Par le th\'eor\`eme de Harrington, Kechris, et Louveau 
(cf [HKL]), on trouve une injection continue $f$ de 
$2^\omega$ dans $X$ telle que $E_0 = (f \times f)^{-1}(A)$ ($E_0$ 
\'etant l'ensemble des suites infinies de 0 et de 1 \'egales \`a 
partir d'un certain rang).\bigskip

 Il suffit donc de regarder le cas o\`u $A = E_0$. Posons, si 
$\alpha \in 2^\omega$, $g_0(\alpha)(p) = 1-\alpha(p)$, et 
$g_n(\alpha)(p) = \alpha(p) \Leftrightarrow p > (n-1)$ si 
$n > 0$ ; la suite $(g_n)$ est une suite d'hom\'eomorphismes, 
et converge simplement vers $g$. Posons \'egalement $B := 
\bigcup_{n\in\omega}  Gr(g_n)$, $C := \check Gr(g_0)$ ; 
par ce qui pr\'ec\`ede $B$ est bor\'elien \`a coupes compactes, 
donc est $\mbox{pot}(\bormone)$ ; $C$ est ouvert, et on a la 
double \'egalit\'e $B \cap A = C \cap B = 
\bigcup_{n\in\omega\setminus\{0\}} Gr(g_n)$. Ce dernier 
ensemble est non-$\mbox{pot}(\bormone)$, \`a cause du lemme 3.5 
appliqu\'e \`a $X = Y = C_n = D_n = 2^\omega$, $(g_n)$, et 
$A = B$.$\hfill\square$\bigskip

 On donne maintenant un r\'esultat qui est une condition 
suffisante pour obtenir la conclusion du (2) de la preuve 
pr\'ec\'edente, mais sans borne sur la classe de Wadge 
potentielle de $A$ (par la remarque 1.5.(a), un bor\'elien \`a 
coupes verticales co-d\'enombrables est $\mbox{pot}(\bormtwo)$).\bigskip

Ce r\'esultat se rattache aussi au th\'eor\`eme 2.3, qui implique 
que si $A$ est non-$\mbox{pot}(\boraone)$, $A$ contient un 
graphe de fonction injective continue d\'efinie sur une 
copie de $2^\omega$ (la r\'eciproque \'etant fausse : prendre 
$A = 2^\omega \times 2^\omega$ ! ). 
Cependant, la r\'eciproque est vraie si la trace de $A$ sur 
un rectangle parfait est un tel graphe, comme on l'a vu 
dans la preuve du (ii) de la proposition 4.8.

\begin{thm} Soient $X$ et $Y$ des espaces 
polonais et $A$ un bor\'elien \`a coupes horizontales maigres 
de $X \times Y$ tel que $\Pi_X''A$ soit non maigre ; alors 
il existe une copie $P$ (resp. $Q$) de $2^\omega$ dans $X$ 
(resp. $Y$), et un hom\'eomorphisme $\phi$ de $P$ sur $Q$ 
tels que $Gr(\phi) = (P \times Q) \cap A$.\end{thm}

 Rappelons un lemme d\'emontr\'e dans dans [Ke] :

\begin{lem} Si $X$ et $Y$ sont des espaces 
polonais et $A$ est bor\'elien \`a coupes horizontales maigres 
de $X \times Y$, $A$ est contenu dans une r\'eunion 
d\'enombrable de bor\'eliens \`a coupes ferm\'ees rares de 
$X \times Y$.\end{lem}

\begin{lem} Soient $G$ et $Y$ des espaces 
polonais, $f$ une fonction continue ``meager-to-one" de 
$G$ dans $Y$, $O$ un ouvert \`a coupes horizontales denses 
de $G \times Y$, $\varepsilon > 0$, et pour $i = 0$, $1$, $M_i$ 
(resp. $N_i$) des ouverts non vides de $G$ (resp. $Y$) 
tels que $M_i \subseteq f^{-1}(N_i)$. Alors il existe des 
ouverts non vides $M'_i$ (de $G$) et $N'_i$ (de $Y$) tels que 
$$\begin{array}{ll}
& (1)~~\overline{M'_i \times N'_i} \subseteq M_i \times N_i \cr 
& (2)~~\delta(M'_i \times N'_i) < \varepsilon \cr 
& (3)~~N'_0 \cap N'_1 = \emptyset \cr 
& (4)~~M'_i \times N'_{1-i} \subseteq O \cr 
& (5)~~M'_i \subseteq f^{-1}(N'_i)
\end{array}$$\end{lem}

\noindent\bf D\'emonstration.\rm\ $f$ est ``meager-to-one", donc on peut 
trouver $y_i$ dans $f''M_i$, $y_0 \not= y_1$, et des 
ouverts $O_i$ de $Y$ avec $y_i \in O_i \subseteq 
\overline{O_i} \subseteq N_i$ et $O_0 \cap O_1 = \emptyset$.\bigskip

 Posons $R_0 := (M_0 \cap f^{-1}(O_0)) \times 
\Pi_Y''[(M_1 \times O_1) \cap Gr(f)]$ ; remarquons que $R_0$ 
est non vide. En effet, si $x_0$ est un ant\'ec\'edent de $y_0$ 
dans $M_0$, $(x_0,y_1)$ est dans $R_0$. On trouve alors 
$(x,y)$ dans $R_0\cap O$ : en effet, $R_0 \setminus O$ est 
ferm\'e de $R_0$ \`a coupes rares relativement \`a $M_0 \cap 
f^{-1}(O_0)$ (qui est ouvert de $G$), donc est ferm\'e rare 
de $R_0$. Soient $L_1$ et $N'_1$ des ouverts de diam\`etre 
au plus $\varepsilon$ tels que $(x,y) \in L_1 \times N'_1 \subseteq 
\overline{L_1} \times N'_1 \subseteq 
[(M_0 \cap f^{-1}(O_0)) \times O_1] \cap O$.

\vfill\eject

  Alors posons $R_1 := (M_1 \cap f^{-1}(N'_1)) \times 
\Pi_Y''[(L_1Ê\times O_0) \cap Gr(f)]$ ; l\`a encore, $R_1$ 
est non vide : en effet, si $x_1$ est un ant\'ec\'edent de $y$ 
dans $M_1$, $(x_1,f(x))$ est dans $R_1$. On trouve $(x',y')$ 
dans $R_1 \cap O$, comme avant, et aussi des ouverts $M'_1$, 
$N'_0$, de diam\`etre au plus $\varepsilon$, tels que l'on ait 
$$(x',y') \in M'_1 \times N'_0 \subseteq 
\overline{M'_1} \times N'_0 \subseteq [(M_1 \cap f^{-1}(N'_1)) 
\times O_0] \cap O.$$ 
Il reste \`a poser $M'_0 := L_1 \cap f^{-1}(N'_0)$.$\hfill\square$\bigskip

\noindent\bf D\'emonstration du th\'eor\`eme 4.10.\rm\ Soit $N$ un ouvert non 
vide de $X$ sur lequel l'analytique $\Pi_X''A$ est 
co-maigre. On peut alors appliquer le lemme 4.11 \`a $N$, 
$Y$, et $A \cap (N \times Y)$, ce qui fournit une suite 
croissante $(F_n)$ de bor\'eliens \`a coupes ferm\'ees rares de 
$N \times Y$, dont la r\'eunion contient $A \cap (N \times Y)$. 
Par le rappel 1.1.(b) on trouve une topologie polonaise $\tau$ 
sur $Y$ affinant la topologie initiale de sorte que $F_n$ 
soit ferm\'e de $N \times (Y, \tau)$.\bigskip

 Par le th\'eor\`eme de Jankov-von Neumann, on trouve une fonction 
$f$ Baire-mesurable uniformisant $A \cap (N \times Y)$ sur 
$N \cap \Pi_X''A$ ; soit alors $G$ un $G_{\delta}$ dense de 
$N$, contenu dans $N \cap \Pi_X''A$, sur lequel $f$ est 
continue ($G$ existe car $N \cap \Pi_X''A$ est co-maigre 
dans $N$).\bigskip

 On construit des ouverts non vides $(V_s)_{s\in 2^{<\omega}}$ 
de $G$ et $(W_s)_{s\in 2^{<\omega}}$ de $(Y, \tau)$ v\'erifiant, si 
$$U_s := V_s \times W_s\mbox{,}$$
$$\begin{array}{ll} 
& (i)~~~\overline{U_{s^{\frown} i}} \subseteq U_s \cr 
& (ii)~~\delta(U_s) < \vert s\vert^{-1}~\mbox{si}~s \not= \emptyset \cr 
& (iii)~W_{s^{\frown} 0} \cap W_{s^{\frown} 1} = \emptyset \cr 
& (iv)~~ \forall~s \not= t \in 2^{n+1}~V_s \times W_t \subseteq \check F_n \cr 
& (v)~~~V_s \subseteq f^{-1}(W_s)
\end{array}$$

 Si on a construit ces objets, posons 
$$\{(x_{\alpha},y_{\alpha})\} := 
\bigcap_{n\in\omega} U_{\alpha \lceil n} = \bigcap_{n\in\omega} 
\overline{U_{\alpha \lceil n}}.$$ 
Alors $\{(x_{\alpha},y_{\alpha})~/~\alpha \in 2^\omega \}$ 
est une copie de $2^\omega$ qui le graphe d'une injection 
partielle $\phi$ (ceci r\'esulte de (iii) et (v), qui assure la 
disjonction de $V_{s^{\frown} 0}$ et $V_{s^{\frown} 1}$ si 
$s$ est dans $2^{<\omega}$). Les projections de ce graphe 
d\'efinissent les copies $P$ et $Q$ de $2^\omega$ annonc\'ees. 
Par (v), $\phi$ est la restriction de $f$ \`a $P$, qui est 
contenu dans $G$, donc $\phi$ est un hom\'eomorphisme de $P$ sur $Q$. 
L'inclusion du graphe de $\phi$ dans $(P \times Q) \cap A$ est 
alors claire. Inversement, si $(x_{\alpha},y_{\beta})$ est 
dans $(P \times Q) \setminus Gr(\phi)$, $\alpha \not= \beta$ 
donc par (iv) on trouve $n > 0$ tel que $\alpha \lceil n \not= 
\beta \lceil n$ ; par cons\'equent $(x_{\alpha},y_{\beta})$ 
n'est pas dans $F_m$ si $n \leq m+1$ et $(x_{\alpha},y_{\beta})$ 
n'est pas dans $A$.\bigskip

  Montrons donc que la construction est possible. On pose $V_{\emptyset} := G$ 
et $W_{\emptyset} := Y$. Admettons avoir construit $U_s$ 
pour $s$ dans $2^{\leq n}$ v\'erifiant les conditions (i)-(v).\bigskip

 On va appliquer plusieurs fois le lemme 4.12, \`a $f$ et \`a 
des ouverts de $G$ et $(Y,\tau)$, ce qui est licite 
puisque $f$ est ``meager-to-one" sur $N$, donc sur $G$ qui 
est co-maigre dans $N$.

\vfill\eject

 On commence \`a appliquer ce lemme 4.12 \`a $\varepsilon = (n+1)^{-1}$, 
$O = (G \times Y) \setminus F_n$, $M_0 = M_1 = V_s$ et 
$N_0 = N_1 = W_s$, de sorte qu'on a assur\'e (i)-(iii) et 
(v), ainsi que (iv) pour les couples de la forme 
$(s^{\frown}i,s^{\frown }(1-i))$, avec $s$ dans $2^n$ et 
$i$ dans 2 (on obtient ainsi $\tilde V_{s^{\frown}0}$, 
$\tilde V_{s^{\frown}1}$, $\tilde W_{s^{\frown}0}$, 
$\tilde W_{s^{\frown}1}$). 
Si jamais, pour $u$ et $v$ distincts dans $2^{n+1}$, 
$(\tilde V_u \times \tilde W_v) \cap F_n$ est non vide, on 
diminue \`a l'aide du lemme 4.12 appliqu\'e \`a $M_0 = \tilde V_u$, 
$M_1 = \tilde V_v$, $N_0 = \tilde W_u$, $N_1 = \tilde W_v$, 
$\tilde V_u$ et $\tilde W_v$ de mani\`ere \`a \'eviter $F_n$ 
tout en conservant (v). On r\'ealise ainsi (iv) au bout 
d'un nombre fini de changements \'eventuels (major\'e par 
$(2^{n+1})^2$).$\hfill\square$\bigskip

 Par le (i) de la proposition 4.8, on a $P_1(\mbox{pot}(\boratwo)\cap \mbox{pot}(\bormtwo))$, 
donc on a fortiori $P_1(\mbox{pot}(\check D_2(\boraone)))$. On va 
donner une nouvelle preuve de ceci mais sous une forme 
beaucoup plus forte, du type du th\'eor\`eme rappel\'e au d\'ebut 
du paragraphe.\bigskip

 Comme il r\'esulte de la preuve du (ii) de la proposition 
4.8, si $D$ est \`a coupes d\'enombrables et est non-$\mbox{pot}(\boraone)$, 
on trouve des injections continues $\phi$ et $\psi$ telles 
que $\Delta (2^\omega) = (\phi \times \psi)^{-1}(D)$, donc 
en particulier $\Delta (2^\omega) \leq_P A$. Mais dans la 
preuve du th\'eor\`eme 3.7, on a vu que $\Delta (2^\omega) 
\not\leq_P \check D_0$, et par 3.6, $\check D_0$ est non-$\mbox{pot}(\boraone)$ ; 
mais pour montrer que 
$\Delta (2^\omega) \not\leq_P \check D_0$, on a utilis\'e le 
fait que $\check D_0$ est \`a coupes co-d\'enombrables. Or il se 
trouve que les $\mbox{pot}(\bormone)$ non-$\mbox{pot}(\boraone)$ \`a coupes 
co-d\'enombrables n'existent pas : sinon soit $\tau$ une 
topologie sur $Y$ rendant les coupes d'un tel $A$ ferm\'ees ; 
$\Pi_Y''(\check A)$ serait ouvert de $(Y, \tau)$ et on 
aurait $\Pi_Y''(\check A) = \bigcup_{x\in X} \check A(x) = 
\bigcup_{n\in\omega} \check A(x_n)$ car $\Pi_Y''(\check A)$ 
est un espace de Lindel{\"o}f ; comme $\check A(x_n)$ est 
d\'enombrable, $\Pi_Y''(\check A)$ le serait aussi et par la 
remarque 2.1, $A$ serait $\mbox{pot}(\borone)$. On peut donc se 
demander si on n'a pas l\`a une caract\'erisation des ``vrais" 
$\mbox{pot}(\bormone)$, \`a savoir : si $B$ est $\mbox{pot}(\bormone)$, $B$ 
est non-$\mbox{pot}(\boraone)$ si et seulement si 
$\Delta (2^\omega) \leq_P B$. On va voir que c'est bien le 
cas. On note 
$L_0 := \{ (\alpha,\beta) \in 2^\omega \times 2^\omega~/~\alpha 
\leq_{\mbox{lex}} \beta \}$.

\begin{thm} Si $A$ est $\mbox{pot}(D_2(\boraone))$ 
dans un produit de deux espaces polonais, les conditions 
suivantes sont \'equivalentes :\smallskip

\noindent (i) $A$ est non-$\mbox{pot}(\boraone)$.\smallskip

\noindent (ii) Il existe des fonctions injectives continues $\phi$, de 
$2^\omega$ dans $X$, et $\psi$, de $2^\omega$ dans $Y$, 
telles que $\Delta (2^\omega) = (\phi \times \psi)^{-1}(A)$ 
ou $L_0 = (\phi \times \psi)^{-1}(A)$.\end{thm}

\noindent\bf D\'emonstration.\rm\ Montrons que (ii) implique (i). On a vu 
au chapitre 2 que $\Delta (2^\omega)$ est non-$\mbox{pot}(\boraone)$, 
donc si $\Delta (2^\omega) = (\phi \times \psi)^{-1}(A)$, $A$ 
est non-$\mbox{pot}(\boraone)$. Il suffit donc de voir que $L_0$ 
est non-$\mbox{pot}(\boraone)$, et il suffit de montrer que 
$\Delta (2^\omega) \leq_P L_0$. Posons donc :
$$f(\alpha,\beta) = \left\{\!\!\!\!\!\!
\begin{array}{ll} 
& (\alpha,\beta)~\mbox{si}~\alpha \geq \beta\mbox{,}\cr 
& (\beta,\alpha)~\mbox{sinon.}
\end{array}\right.$$
Alors $\Delta (2^\omega) = f^{-1}(L_0)$, et $f$ est dans $C_0$ :
$$\begin{array}{ll} 
f(\alpha,\beta) \in C \times D 
& \Leftrightarrow~(\alpha \geq \beta~\mbox{et}~\alpha \in C~\mbox{et}~\beta \in D)
~\mbox{ou}~(\alpha < \beta~\mbox{et}~\alpha \in D~\mbox{et}~\beta \in C)\cr 
& \Leftrightarrow~(\alpha \in C \cap D~\mbox{et}~\beta \in C \cap D)
~\mbox{ou}~(\alpha > \beta~\mbox{et}~\alpha \in C~\mbox{et}~\beta \in D)~\mbox{ou}~\cr 
& ~~~~~~(\alpha < \beta~\mbox{et}~\alpha \in D~\mbox{et}~\beta \in C)
\end{array}$$
 ~~~~Montrons maintenant que (i) implique (ii) : $A$ est 
$\mbox{pot}(D_2(\boraone))$, donc par la proposition 2.2 on a 
$A = (\bigcup_{n\in\omega} C_n \times D_n) \setminus V$, 
avec $C_n$, $D_n$ bor\'eliens, $V$ dans $\mbox{pot}(\boraone)$, et 
la r\'eunion disjointe.

\vfill\eject

 En effet, on peut supposer cette r\'eunion des rectangles 
disjointe : si $E = \bigcup_{n\in\omega}  A_n \times B_n$,  
et si $I \subseteq n$, posons : $C_{n,I} := A_n \cap 
\bigcap_{i\in I} [A_i \setminus (\bigcup_{j\in(n\setminus I)} A_j)]$, 
$D_{n,I} := B_n \setminus (\bigcup_{i\in I} B_i)$. Alors 
on v\'erifie facilement que $E = \bigcup_{n,I\subseteq n} C_{n,I} 
\times  D_{n,I}$, cette r\'eunion \'etant disjointe.\bigskip

 On trouve alors $n$ tel que $(C_n \times D_n) \setminus V$ 
ne soit pas $\mbox{pot}(\boraone)$ ; or $(C_n \times D_n) \setminus V$ 
est $\mbox{pot}(\bormone)$, donc si on admet le r\'esultat pour $A$ 
dans $\mbox{pot}(\bormone)$, on a le th\'eor\`eme, par disjonction de 
la r\'eunion (en effet, $\phi$ (resp. $\psi$) est \`a valeurs 
dans $C_n$ (resp. $D_n$)).\bigskip

 Quitte \`a affiner les topologies, on peut supposer que $A$ 
est ferm\'e, puisque la continuit\'e de $\phi$ et $\psi$ avec 
des topologies plus fines entrainera la continuit\'e avec 
les topologies initiales.\bigskip

 On construit par r\'ecurrence sur $\vert s\vert$, o\`u $s\in 2^{<\omega}$, 
des ouverts-ferm\'es non vides $V_s \subseteq X$, $W_s \subseteq 
Y$ tels que si $U_s := V_s \times W_s$,
$$\begin{array}{ll} 
& (i)~~~U_{s^\frown i} \subseteq U_s \cr 
& (ii)~~\delta(U_s) < \vert s\vert^{-1}~\mbox{si}~s \not= \emptyset \cr 
& (iii)~V_{s^\frown 0} \cap V_{s^\frown 1} = W_{s^\frown 0} \cap W_{s^\frown 1} = \emptyset \cr 
& (iv)~~U_s \cap A \notin \mbox{pot}(\boraone) \cr 
& (v)~~~V_{s^\frown 0} \times W_{s^\frown 1} \subseteq \check A
\end{array}$$
Posons $U_\emptyset := X \times Y$, et admettons avoir 
construit $(U_s)_{s\in 2^{\leq n}}$ v\'erifiant (i)-(v). 
Recouvrons $V_s$ et $W_s$ par des ouverts-ferm\'es de 
diam\`etre au plus $(\vert s\vert +1)^{-1}$, soient $(V_n)$ et 
$(W_p)$. On a 
$$A \cap U_s = \bigcup_{n,p\in\omega} A \cap 
(V_n \times W_p)\mbox{,}$$ 
donc l'un des $A \times (V_n \times W_p)$ 
est non-$\mbox{pot}(\boraone)$. Posons $V := V_n$, $W := W_p$ ; 
on a $U := V \times W \subseteq U_s$, $A\cap U \notin 
\mbox{pot}(\boraone)$, et $\delta(U) < (\vert s\vert +1)^{-1}$.\bigskip 

 Posons $C := \cup ~\{B~/~B \in \boraone \lceil U$ et $A \cap 
B \in \mbox{pot}(\boraone)\}$, $A' := (A \cap U) \setminus C$. 
Comme $C$ est de Lindel\"of, $A \cap C$ est $\mbox{pot}(\boraone)$.\bigskip

 Montrons qu'il existe $(x,x',y,y')$ dans $V^2 \times W^2$ 
tels que $(x, y)$ et $(x',y')$ soient dans $A'$ et $(x', y)$ 
ne soit pas dans $A$. Si tel n'est pas le cas, 
$\Pi_X''A' \times \Pi_Y''A' \subseteq U \cap A$, donc comme 
$V$, $W$, et $A$ sont ferm\'es, 
$\overline{\Pi_X"A'} \times \overline{\Pi_Y"A'} \subseteq 
U \cap A$, d'o\`u l'\'egalit\'e 
$$A \cap U = [\overline{\Pi_X''A'} \times \overline{\Pi_Y''A'}] 
\cup [A \cap ((V \setminus \overline{\Pi_X''A'}) \times  W)] \cup 
[A \cap (V \times (W \setminus \overline{\Pi_Y''A'}))].$$
Le premier terme du membre de droite est un rectangle 
bor\'elien, donc est $\mbox{pot}(\boraone)$ ; le deuxi\`eme est 
contenu dans $C$, et vaut $A \cap C \cap [(V \setminus \overline{\Pi_X''A'}) \times  W]$, 
donc est $\mbox{pot}(\boraone)$, ainsi que le troisi\`eme. Donc 
$A \cap U$ est $\mbox{pot}(\boraone)$, ce qui est exclus.\bigskip

 Comme $U \setminus A$ est ouvert, soient $Z$, $T$ des 
ouverts tels que $(x',y) \in Z \times T \subseteq U \setminus A$ ; 
on peut poser $V_{s^\frown 0} := Z$, $V_{s^\frown 1} := 
V \setminus Z$, $W_{s^\frown 0} := W \setminus T$, $W_{s^\frown 1} := T$.\bigskip

 V\'erifions (iv) : $(x',y')$ est dans $U_{s^\frown 0} \cap A 
= (Z \times W) \cap A$, donc $U_{s^\frown 0} \cap A$ n'est 
pas $\mbox{pot}(\boraone)$. De m\^eme pour $U_{s^\frown 1} \cap A$, 
avec $(x,y)$.\bigskip

 Soit $\Phi$, de $2^\omega$ dans $X \times Y$, d\'efinie par : 
$\{ \Phi(\alpha)\} = \bigcap_{n\in\omega} U_{\alpha \lceil n}$ ; 
alors $\Phi''2^\omega = Gr(f)$, o\`u $f$ est une injection 
continue d\'efinie sur la copie $\Pi_X''Gr(f)$ de $2^\omega$. 
De plus $Gr(f) \subseteq A$ car par (iv), il existe 
$(\gamma_n^\alpha,\beta_n^\alpha)$ dans $U_{\alpha \lceil n} \cap A$, 
et la suite 
$((\gamma_n^\alpha,\beta_n^\alpha))_{n\in\omega}$ converge 
vers $\Phi(\alpha)$, donc $\Phi(\alpha)$ est dans $A$ qui 
est ferm\'e. Enfin, si $\alpha <_{\mbox{lex}} \beta$, 
$((\Pi_X \circ \Phi)(\alpha),(\Pi_Y \circ \Phi)(\beta))$ 
n'est pas dans $A$, par (v). Posons donc 
$E := ((\Pi_X\circ \Phi) \times (\Pi_Y\circ \Phi))^{-1}(A)$ ; on a 
$$\begin{array}{ll} 
& (a)~\Delta(2^\omega) \subseteq E \cr 
& (b)~\alpha <_{\mbox{lex}} \beta \Rightarrow (\alpha,\beta) \notin E
\end{array}$$
 Montrons que $2^\omega$ contient une copie $Z$ de $2^\omega$ 
telle que $E \cap Z^2 = \{ (\alpha ,\beta )\in 2^\omega\times 2^\omega /(\beta ,\alpha )\in L_0\}\cap Z^2$ ou 
$\Delta(2^\omega) \cap Z^2$.\bigskip

\noindent\bf Premier cas.\rm\ Pour toute suite de $2^{<\omega}$, 
$N_s^2 \setminus E$ n'est pas antisym\'etrique.\bigskip

 On construit alors par r\'ecurrence sur $\vert s\vert$ une suite 
$(V_s)$ de $\borone (2^\omega )$ non vides v\'erifiant 
$$\begin {array}{ll} 
& (1)~V_{s^\frown i} \subseteq V_s \cr 
& (2)~\delta(V_s) < \vert s\vert^{-1}~\mbox{si}~s\not= \emptyset \cr 
& (3)~(V_{s^\frown 0} \times V_{s^\frown 1}) 
\cup (V_{s^\frown 1} \times V_{s^\frown 0}) \subseteq \check E
\end{array}$$
Ceci ne pose aucun probl\`eme puisque $E$ est ferm\'e dans 
$2^\omega\times2^\omega$ et que $V_s^2 \setminus E$ n'est 
pas antisym\'etrique. La formule $\{g(\alpha)\} := \bigcap_{n\in\omega} 
V_{\alpha \lceil n}$ d\'efinit une injection continue $g$, 
qui est un hom\'eomorphisme de $2^\omega$ sur son image $Z$, 
qui par (a) v\'erifie $E \cap Z^2 = \Delta(2^\omega) \cap Z^2$.\bigskip

\noindent\bf Second cas.\rm\ Il existe $s$ dans $2^{<\omega}$ 
telle que $N_s^2 \setminus E$ soit antisym\'etrique.\bigskip

 Alors par (a) et (b), $Z := N_s$ v\'erifie 
$E \cap Z^2 = \{ (\alpha ,\beta )\in 2^\omega\times 2^\omega /(\beta ,\alpha )\in L_0\}\cap Z^2$.\bigskip

\noindent Soit alors $f$ un hom\'eomorphisme d\'ecroissant de $2^\omega$ 
sur $Z$ ; les fonctions compos\'ees 
$$\phi := \Pi_X \circ \Phi \circ f$$ 
et $\psi := \Pi_Y \circ \Phi \circ f$ r\'epondent au 
probl\`eme.$\hfill\square$\bigskip

\begin{cor} (a) Sous les hypoth\`eses du 
th\'eor\`eme 4.13, on a 
$$A \notin \mbox{pot}(\boraone)\Leftrightarrow 
\Delta(2^\omega) \leq_P A.$$
(b) On ne peut pas trouver $A_0$ tel que si $A$ est 
bor\'elien d'un produit de deux espaces polonais, on ait 
$A\notin \mbox{pot}(\boraone)\Leftrightarrow A_0 \leq_P A$.\smallskip

\noindent (c) Soit $A$ un bor\'elien d'un produit de deux espaces 
polonais. Alors $A\in \mbox{pot}(\borone)$ si et seulement si $A <_P \Delta(2^\omega)$.\end{cor}

\noindent\bf D\'emonstration.\rm\ (a) Si $A$ est non-$\mbox{pot}(\boraone)$, par 
le th\'eor\`eme 4.13, $\Delta(2^\omega) \leq_P A$ ou $L_0 \leq_P A$, 
et on a vu dans la preuve que $\Delta(2^\omega) \leq_P L_0$, 
donc $\Delta(2^\omega) \leq_p A$. La r\'eciproque r\'esulte du 
chapitre 2.

\vfill\eject

\noindent (b) Avec $A := \Delta(2^\omega)$, on voit que si 
$A_0$ existait, $A_0$ serait $\mbox{pot}(\bormone)$. Avec $A := A_0$, 
on voit que $A_0$ serait $\mbox{pot}(\bormone)$ non-$\mbox{pot}(\boraone)$. 
Par (a), on aurait donc $\Delta(2^\omega) \leq_p A_0$. 
Avec $A := \check D_0$, on aurait $A_0 \leq_P \check D_0$, donc 
$\Delta(2^\omega) \leq_P \check D_0$, ce qui est exclus par la preuve 
du th\'eor\`eme 3.7.\bigskip

\noindent (c) Si $A$ est $\mbox{pot}(\borone)$, et si 
$$f(x,y) := \left\{\!\!\!\!\!\!
\begin{array}{ll} 
& (0^\omega,0^\omega)~\mbox{si}~ (x,y) \in A\mbox{,}\cr
& (0^\omega,1^\omega)~\mbox{sinon,}
\end{array}\right.$$ 
$f$ est dans $C_0$ et r\'eduit $A$ \`a $\Delta(2^\omega)$ ; $\Delta(2^\omega) \not\leq_P A$ car 
$\Delta(2^\omega)$ n'est pas $\mbox{pot}(\boraone)$.\bigskip

 R\'eciproquement, $A \leq_P \Delta(2^\omega)$, donc $A$ est 
$\mbox{pot}(\bormone)$ ; si $A$ \'etait non-$\mbox{pot}(\boraone)$, par 
(a) on aurait $\Delta(2^\omega) \leq_P A$, ce qui est exclus.$\hfill\square$\bigskip

 Le contre-exemple $\check D_0$ prouve que l'hypoth\`ese ``$A \in 
\mbox{pot}(D_2(\boraone))$" du th\'eor\`eme 4.13 et du corollaire 
4.14.(a) est optimale du point de vue classe de Wadge.\bigskip

 On ne peut pas se ramener \`a un seul exemple typique dans 
le th\'eor\`eme 4.13, avec des r\'eductions rectangulaires : si 
$$C \leq_R D\ \Leftrightarrow\ \mbox{il existe}~f, g\ \mbox{bor\'eliennes telles que}~
C = (f \times g)^{-1}(D)\mbox{,}$$
alors $\Delta(2^\omega) \perp_R L_0$, comme on le v\'erifie imm\'ediatement.

\section{$\!\!\!\!\!\!$ Uniformisation partielle des $G_\delta$.}\indent

 On va montrer dans ce paragraphe des th\'eor\`emes
d'uniformisation partielle des $G_{\delta}$, dont le
but est d'essayer de trouver une r\'eciproque au
lemme 3.5. Ces th\'eor\`emes sont \`a rapprocher d'une part de r\'esultats dans [GM] et [Ma], o\`u
 au lieu de consid\'erer la cat\'egorie, il est question
d'ensembles de mesure 1 sur chacun des facteurs ; et d'autre 
part de r\'esultats de G. Debs et J. Saint Raymond, o\`u il est 
question de fonctions totales et injectives, avec des 
hypoth\`eses de compacit\'e sur chacun des facteurs (cf 
[D-SR]).

\begin{lem} Soient $X'$ un ouvert-ferm\'e non vide 
de $\omega^\omega$, $Y$ un espace polonais, $Y'$ un ouvert non 
vide de $Y$, $\varepsilon > 0$, et $O$ un ouvert dense de 
$X' \times Y'$ dont la projection est $X'$ ; il existe des 
suites $(U_k)$ (d'ouverts-ferm\'es non vides de $\omega^\omega$) et 
$(V_k)$ (d'ouverts non vides de $Y$) telles que 
$$\begin{array}{ll} 
& (1)~\bigcup_{k\in\omega} U_k~(\mbox{resp.}~\bigcup_{k\in\omega} \overline{V_k})~\mbox{est~dense~dans}~X' ~(\mbox{resp.}~Y')\cr 
& (2)~U_k \times V_k \subseteq O\cr 
& (3)~\delta(U_k),~ \delta(V_k) < \varepsilon\cr 
& (4)~U_p \cap U_q = \emptyset~\mbox{si}~p\not= q
\end{array}$$\end{lem}

\noindent\bf D\'emonstration.\rm\ Soit $(W_m)$ une partition de $X'$ en 
ouverts-ferm\'es, avec $\delta (W_m) < \varepsilon$ et $W_m \not= \emptyset$ 
(c'est possible car $X'$ est hom\'eomorphe \`a $\omega^\omega$). Soit 
$(T_m)$ une base de la topologie de $Y'$. Par densit\'e de 
$O$, on trouve $(x_m,y_m)$ dans $(W_m \times T_m) \cap O$, 
et un ouvert-ferm\'e $X_m$ de $X'$, un ouvert $Y_m$ de $Y'$ 
tels que $\delta (Y_m) < \varepsilon$, et $(x_m,y_m) \in X_m 
\times \overline{Y_m} \subseteq (W_m \times T_m) \cap O$. 
Si $\bigcup_{m\in\omega} X_m$ est dense dans $X'$, on a 
construit, en prenant $U_k := X_k$ et $V_k := Y_k$, les 
ouverts cherch\'es ; en effet, ils v\'erifient bien les 
conditions $(1)\mbox{-}(4)$.

\vfill\eject

 Sinon, si $x$ est dans $X'\setminus 
\overline{\bigcup_{m\in\omega} X_m}$, on trouve un 
ouvert-ferm\'e $Z_x$ de $X'$, $y_x$ dans $O(x)$, et un 
ouvert $R_x$ de $Y$ tels que l'on ait $\delta (Z_x)$,  
$\delta (R_x) < \varepsilon$, et \'egalement l'inclusion 
$$(x,y_x) \in Z_x \times \overline{R_x} \subseteq O \cap 
[(X' \setminus \overline{\bigcup_{m\in\omega} X_m}) \times Y'].$$
On a $X' \setminus \overline{\bigcup_{m\in\omega} X_m} = \cup Z_x 
= \bigcup_{n\in\omega} Z_{x_n} = 
\bigcup_{n\in\omega} [ Z_{x_n}\setminus (\bigcup_{p<n} Z_{x_p})]$, 
puisque $X'\setminus \overline{\bigcup_{m\in\omega} X_m}$ 
est un espace de Lindel\"of. Posons $\tilde Z_n := Z_{x_n} \setminus 
(\bigcup_{p<n} Z_{x_p})$ ; s'il y a une infinit\'e de 
$\tilde Z_n$ non vides, on note $(Z_n)$ la suite form\'ee de 
ces ouverts-ferm\'es (on trouve donc $n_k$ tel que $Z_k = \tilde Z_{n_k}$).\bigskip

 Sinon, on partitionne un $\tilde Z_{n_0}$ non vide en une 
suite d'ouverts-ferm\'es non vides de $X'$, et on note $(Z_n)$ 
la suite form\'ee des $\tilde Z_n$ non vides (pour $n \not= n_0$) 
et de la partition de $\tilde Z_{n_0}$ ; on a encore 
$Z_k \subseteq \tilde Z_{n_k}$, avec \'eventuellement 
$n_k = n_{k'}$.\bigskip

  Il reste \`a poser : $R_k = R_{x_{n_k}}$, puis $U_{2k} := X_k$, 
$U_{2k+1} := Z_k$, $V_{2k} := Y_k$, $V_{2k+1} := R_k$.$\hfill\square$

\begin{thm} Soient $X$, $Y$ des espaces 
polonais non vides, $X$ \'etant parfait, $A$ un $G_{\delta}$ 
dense de $X \times Y$. Alors il existe des $G_{\delta}$ 
denses $F$ (dans $X$) et $G$ (dans $Y$) et une surjection 
continue ouverte $f$ de $F$ sur $G$ dont le graphe est 
contenu dans $A$ (avec de plus $F$ hom\'eomorphe \`a $\omega^\omega$).\end{thm}

\noindent\bf D\'emonstration.\rm\ On peut supposer $A$ \`a coupes verticales 
co-maigres : $\check A$ est maigre, donc a ses coupes verticales 
maigres, sauf sur un ensemble maigre (th\'eor\`eme de 
Kuratowski-Ulam). $A$ a donc ses coupes verticales 
co-maigres sur un $G_{\delta}$ dense $H$ de $X$, qui comme 
$X$ est polonais parfait.\bigskip

 On peut supposer \'egalement que $X = \omega^\omega$ (de sorte que 
si $F$ est $G_{\delta}$ dense de $X$, $F$ est hom\'eomor-phe 
\`a $\omega^\omega$) : en effet, on proc\`ede comme dans la preuve du 
corollaire 1.7.\bigskip

 Soit $(O_n)$ une suite d'ouverts denses de $X \times Y$ 
telle que $A = \bigcap_{n\in\omega} O_n$, $\phi_0$ de 
$\omega$ dans $\{\emptyset \}$ et, si $n > 0$, $\phi_n$ 
une bijection de $\omega$ sur $\omega^n$. On construit une 
suite $(U_s)_{s\in \omega^{<\omega}}$ d'ouverts-ferm\'es non 
vides de $X$, et une suite $(V_s)_{s\in \omega^{<\omega}}$ 
d'ouverts non vides de $Y$ v\'erifiant 
$$\begin{array}{ll} 
& (i)~~~\bigcup_{n\in\omega} U_{s^\frown n}~(\mbox{resp.}~
\bigcup_{n\in\omega} \overline{V_{s^\frown n}})~\mbox{est~
dense~dans}~U_s~(\mbox{resp.}~V_s) \cr 
& (ii)~~U_s \times V_s \subseteq 0_{(\vert s\vert -1)}~\mbox{si}~s \not= \emptyset \cr 
& (iii)~\delta(U_s),~\delta(V_s) < \vert s\vert^{-1}~\mbox{si}~s \not= \emptyset \cr 
& (iv)~~U_{s^\frown n} \cap U_{s^\frown m} = \emptyset~\mbox{si}~n\not= m \cr 
& (v)~~~(\bigcup_{k\in\omega} V_{\phi_n(k)}) \cap \bigcup_{m+p<n} 
[V_{\phi_m(p)}\setminus (\bigcup_{l\in\omega} 
V_{\phi_m(p)^\frown l})] = \emptyset
\end{array}$$
On pose $U_\emptyset := \omega^\omega$, $V_\emptyset := Y$ ; 
admettons avoir construit $(U_s)_{\vert s\vert\leq n}$ et 
$(V_s)_{\vert s\vert \leq n}$, 
$(U_{\phi_n(p)^\frown k})_{p<m,k\in\omega}$, 
$(V_{\phi_n(p)^\frown k})_{p<m,k\in\omega}$ v\'erifiant (i)-(v).\bigskip

 \noindent~~~On construit, si ce n'est d\'ej\`a fait, $(U_{\phi_n(m)^\frown k})_{k\in\omega}$, 
$(V_{\phi_n(m)^\frown k})_{k\in\omega}$ en appliquant le 
lemme pr\'ec\'edent \`a  $\varepsilon := (n+1)^{-1}$, 
$X' := U_{\phi_n(m)}$, 
$$Y' := V_{\phi_n(m)} \setminus 
\overline{\bigcup_{q+p<n} [V_{\phi_q(p)}\setminus (\bigcup_{l\in\omega} V_{\phi_q(p)^\frown l})]}
\mbox{,}$$ 
$O := O_n \cap (X'\times Y')$ ; les conditions demand\'ees 
sont v\'erifi\'ees, la densit\'e ne posant pas de probl\`eme car 
l'adh\'erence ci-dessus est rare.

\vfill\eject

 Posons $F := \bigcap_{n\in\omega} (\bigcup_{s\in\omega^n} U_s)$, 
$G := \bigcap_{n\in\omega} (\bigcup_{s\in\omega^n} V_s)$; 
$F$ et $G$ sont clairement des $G_{\delta}$ denses de $X$ 
et $Y$.\bigskip

  Si $x$ est dans $F$, on trouve pour tout $n$ une unique 
suite $s_n$ de $\omega^n$ telle que $x \in U_{s_n}$ et 
aussi $s_n \prec s_{n+1}$ ; alors $(\overline{V_{s_n}})$ 
d\'efinit $f(x)$ dans $G$, et $f$ est clairement continue 
car $f''(F \cap U_s) \subseteq G \cap V_s$. Montrons 
l'inclusion inverse, ce qui ach\`evera la preuve.\bigskip

 Soit donc $y$ dans $G \cap V_s$, et $p$ tel que 
$s = \phi_{\vert s\vert}(p)$ ; alors il existe un entier $\alpha (\vert s\vert )$ 
tel que $y$ soit dans $V_{s^\frown \alpha(\vert s\vert )}$, sinon 
$y \notin G$ car par (v), on aurait alors $y \notin 
\bigcup_{k\in\omega} V_{\phi_{\vert s\vert +p+1}(k)}$ ; on 
construit comme ceci par r\'ecurrence $\alpha$ dans $N_s$ 
tel que $\{ y \} = \bigcap_{n\in\omega} V_{\alpha \lceil n}$ ; 
mais alors $(U_{\alpha \lceil n})$ d\'efinit $x$ dans $F \cap U_s$ 
tel que $f(x) = y$.$\hfill\square$\bigskip

 On ne peut pas supprimer compl\`etement une des hypoth\`eses, 
ni esp\'erer mieux avec ces hypoth\`eses, dans le sens suivant :\bigskip

\noindent - Si on ne suppose pas $X$ parfait, prendre $X = \omega$ et 
$Y = 2^\omega$.\bigskip

\noindent - Si on ne suppose pas que $A$ est $G_{\delta}$, prendre $X = Y = 2^\omega$ 
et $A = (2^\omega \times 2^\omega) \setminus 
(P_{\infty} \times P_{\infty})$.\bigskip

\noindent - Si on ne suppose pas $A$ dense, prendre $X = Y = 2^\omega$ 
et $A = \{(0^\omega, 0^\omega)\}$.\bigskip

\noindent - On ne peut pas avoir $f$ totale ou surjective sur $Y$ : 
prendre $X = Y = 2^\omega$ et $A = P_{\infty}^2$.\bigskip

\noindent - Enfin, on ne peut pas avoir $f$ injective : prendre 
$X = 2^\omega$ et $Y = \omega$.

\begin{lem} Soient $\varepsilon > 0$, $U$ et $V$ des 
ouverts non vides de $\omega^\omega$, et $O$ un ouvert dense de 
$U \times V$ ; on trouve alors des suites $(Z_n)$ et $(T_n)$ 
d'ouverts-ferm\'es non vides de $\omega^\omega$ v\'erifiant 
$$\begin{array}{ll} 
& (i)~~~\delta (Z_n),~\delta (T_n) < \varepsilon \cr 
& (ii)~~Z_n \times T_n \subseteq O \cr 
& (iii)~Z_n \cap Z_m = T_n \cap T_m = \emptyset~\mbox{si}~n\not= m \cr 
& (iv)~~\bigcup_{n\in\omega} Z_n~\mbox{est~dense~dans } ~U
\end{array}$$\end{lem}

\noindent\bf D\'emonstration.\rm\ Soient $(U_n)$ une base de la topologie 
de $U$, et $(V_n)$ une partition de $V$ en ouverts-ferm\'es non 
vides de diam\`etre au plus $\varepsilon$. Alors $O \cap (U_n \times V_n)$ 
est non vide, donc contient $(x_n,y_n)$ et on trouve des 
suites $(X_n)$ et $(Y_n)$ d'ouverts-ferm\'es de $\omega^\omega$, 
avec $\delta (X_n) < \varepsilon$ et 
$$(x_n,y_n) \in X_n \times Y_n 
\subseteq O \cap (U_n \times V_n).$$ 
R\'eduisons la suite $(X_n)$, ce qui fournit $(W_n)$. S'il y a une infinit\'e de 
$W_n$ non vides, c'est termin\'e. Sinon on partitionne un 
$W_{n_0}$ non vide et le $Y_{n_0}$ correspondant en une 
suite infinie d'ouverts-ferm\'es non vides.$\hfill\square$

\begin{thm} Sous les hypoth\`eses du th\'eor\`eme 5.2, si de plus $Y$ est parfait, on peut avoir pour $f$ 
un hom\'eomorphisme.\end{thm}

\vfill\eject

\noindent\bf D\'emonstration.\rm\ Comme dans la preuve du th\'eor\`eme 5.2, 
on peut supposer que $X$, $Y = \omega^\omega$.\bigskip

 Soit $(O_n)$ une suite d'ouverts denses de 
$\omega^\omega \times \omega^\omega$ telle que $A = \bigcap_{n\in\omega} O_n$. 
On construit alors des suites d'ouverts-ferm\'es non vides 
de $\omega^\omega$, $(Z_s)_{s\in \omega^{<\omega}}$ et 
$(T_s)_{s\in \omega^{<\omega}}$ v\'erifiant 
$$\begin{array}{ll} 
& (1)~\delta (Z_s),~\delta (T_s) < 
\vert s\vert^{-1}~\mbox{si}~s \not= \emptyset \cr 
& (2)~Z_s \times T_s \subseteq O_{\vert s\vert -1}~\mbox{si}~s \not= \emptyset \cr 
&   (3)~\bigcup_{n\in\omega} Z_{s^\frown n}~
(\mbox{resp.}~\bigcup_{n\in\omega} T_{s^\frown n})~\mbox{est~dense~dans}~Z_s~
(\mbox{resp.}~T_s) \cr 
& (4)~Z_{s^\frown n} \cap Z_{s^\frown m} = T_{s^\frown n} \cap 
T_{s^\frown m} = \emptyset~\mbox{si}~n \not= m
\end{array}$$
On pose pour commencer $Z_\emptyset = T_\emptyset = \omega^\omega$. 
Admettons avoir construit $(Z_s)_{\vert s\vert\leq n}$ et 
$(T_s)_{\vert s\vert \leq n}$ v\'erifiant $(1)\mbox{-}(4)$.\bigskip

 On applique le lemme 5.3 \`a $\varepsilon := (n+1)^{-1}$, 
$U := Z_s$, $V := T_s$, $O := O_n \cap (U \times V)$. 
Deux cas se pr\'esentent : si $T := \bigcup_{m\in\omega} T_m$ 
est dense dans $T_s$, c'est termin\'e : on pose $Z_{s^\frown m} 
:= Z_m$ et $T_{s^\frown m} := T_m$.\bigskip

 Sinon, on applique le lemme 5.3 \`a 
$O := \{(x,y)~/~(y,x) \in O_n \cap (Z_s \times (T_s\setminus \overline{T}))\}$, 
$U := T_s \setminus \overline{T}$, $V := Z_s$, 
$\varepsilon := (n+1)^{-1}$, ce qui fournit $(Z'_m)$ et $(T'_m)$ 
v\'erifiant 
$$\begin{array}{ll} 
& (i)~~~\delta (Z'_m),~\delta (T'_m) < (n+1)^{-1} \cr 
& (ii)~~Z'_m \times T'_m \subseteq O_n \cap (Z_s \times (T_s\setminus \overline{T})) \cr 
& (iii)~\bigcup_{m\in\omega} T'_m~\mbox{est~dense~dans}~T_s\setminus \overline{T} \cr 
& (iv)~~Z'_p \cap Z'_m = T'_p \cap T'_m = \emptyset~\mbox{si}~p \not= m
\end{array}$$
Posons, si $j \in\omega$, 
$n_j := \mbox{min}~\{m \in\omega~/~Z'_j \cap Z_m \not= \emptyset \}$ ; 
si $m \in\omega$, $I_m := \{j \in\omega~/~n_j = m\}$ ; et 
$L_m := \overline{\bigcup_{j\in I_m} Z'_j \cap Z_m}$. 
Quatre cas se pr\'esentent.\bigskip

\noindent\bf Premier cas.\rm\ $I_m = \emptyset$.\bigskip

 On pose $Z_0^m := Z_m$,~$T_0^m := T_m$, $Z_j^m = T_j^m := \emptyset~\mbox{si}~j\geq 1$.\bigskip

\noindent\bf Deuxi\`eme cas.\rm\ $I_m \not= \emptyset$ et $Z_m = L_m$.\bigskip

 Soit $j_m :=\mbox{min}~I_m$. On partitionne $Z'_{j_m} \cap Z_m$ en 2 ouverts-ferm\'es non vides 
$\tilde Z_0$ et $\tilde Z_1$, et on pose 
$$\begin{array}{ll} 
& Z_0^m := \tilde Z_0,~T_0^m := T_m\mbox{,}\cr 
& Z_1^m := \tilde Z_1,~T_1^m := T'_{j_m}\mbox{,}\cr 
& Z_{j+1}^m := Z_m \cap Z'_j,~T_{j+1}^m := T'_j~\mbox{si}~j > j_m~
\mbox{et}~j \in I_m\mbox{,}\cr 
& Z_j^m = T_j^m := \emptyset~\mbox{si}~j > 1~\mbox{et}~
(j-1 = j_m~\mbox{ou}~j-1 \notin I_m).
\end{array}$$ 
\bf Troisi\`eme cas.\rm\ $I_m \not= \emptyset$ est fini et $L_m \subset_{\not=} Z_m$.\bigskip
 
 On pose 
$$\begin{array}{ll} 
& Z_0^m := Z_m\setminus L_m,~T_0^m := T_m\mbox{,}\cr 
& Z_{j+1}^m := Z_m \cap Z'_j,~T_{j+1}^m := T'_j~\mbox{si}~j \in I_m\mbox{,}\cr 
& Z_j^m = T_j^m := \emptyset~\mbox{si}~j > 1~\mbox{et}~j-1 \notin I_m.
\end{array}$$ 
\bf Quatri\`eme cas.\rm\ $I_m$ est infini et $L_m \subset_{\not=} Z_m$.\bigskip 

 Soient $\psi$ une bijection de $\omega$ sur $I_m$, et 
$(\tilde Z_m)$, $(\tilde T_m)$ des partitions en 
ouverts-ferm\'es non vides de $Z_m \setminus L_m$ et $T_m$. On pose : 
$Z_{2k}^m := \tilde Z_{2k}$,~$T_{2k}^m := \tilde T_k$, $Z_{2k+1}^m := Z_m \cap Z'_{\psi(k)}$,~
$T_{2k+1}^m := T'_{\psi(k)}$.\bigskip 

 On renum\'erote, de fa\c con \`a ce que $\{Z_{s^\frown p}~/~p \in \omega \} 
= \{Z^i_j~/~i,~j\in\omega~{\rm et}~Z^i_j \not= \emptyset \}$ ; 
on pose $T_{s^\frown p} := Z^i_j$ si $Z_{s^\frown p} = Z^i_j$ 
et $(Z_{s^\frown p})$, $(T_{s^\frown p})$ r\'epondent au 
probl\`eme, comme on le v\'erifie facilement.\bigskip

 Il est maintenant clair que l'ensemble 
$\bigcap_{n\in\omega} (\bigcup_{s\in \omega^n} (Z_s\times T_s))$ 
est le graphe d'un hom\'eomorphis-me de $F := \bigcap_{n\in\omega} (\bigcup_{s\in \omega^n} Z_s)$ sur 
$G := \bigcap_{n\in\omega} (\bigcup_{s\in \omega^n} T_s)$.$\hfill\square$\bigskip

 Pour pouvoir appliquer ces r\'esultats aux classes de Wadge 
potentielles, il faut traiter le cas o\`u le $G_{\delta}$ 
est rare.

\section{$\!\!\!\!\!\!$ R\'ef\'erences.}

\noindent [D-SR]\ \ G. Debs et J. Saint Raymond, \it S\'elections
bor\'eliennes injectives,\rm ~Amer. J. Math.~111 (1989), 519-534

\noindent [GM]\ \ S. Graf and R. D. Mauldin,~\it Measurable one-to-one
selections and transition kernels,\rm ~Amer. J. Math.~107 (1985), 407-425

\noindent [HKL]\ \ L. A. Harrington, A. S. Kechris et A. Louveau,~\it A Glimm-Effros 
dichotomy for Borel equivalence relations,\rm~J. Amer. Math. Soc.~3 (1990), 903-928

\noindent [Ke]\ \ A. S. Kechris,~\it Measure and category in effective 
descriptive set theory,~\rm Ann. of Math. logic,~5 (1973), 337-384

\noindent [Ku]\ \  K. Kuratowski,~\it Topology,~\rm Vol. 1, Academic Press, 1966

\noindent [Lo1]\ \ A. Louveau,~\it Livre \`a para\^\i tre\rm

\noindent [Lo2]\ \ A. Louveau,~\it Some results in the Wadge hierarchy of Borel sets,~\rm Cabal Sem. 79-81 
(A. S. Kechris, D. A. Mauldin, Y. N. Moschovakis, eds) Lect. Notes in Math., Springer-Verlag~1019 (1983), 28-55

\noindent [Lo3]\ \ A. Louveau,~\it A separation theorem for $\Ana$ sets,~\rm Trans. A. M. S.~260 (1980), 363-378

\noindent [Lo4]\ \ A. Louveau,~\it Ensembles analytiques et bor\'eliens dans les espaces produit,~\rm 
Ast\'erisque (S. M. F.)~78 (1980)

\noindent [Lo-SR]\ \ A. Louveau and J. Saint Raymond,~\it Borel classes and closed games : Wadge-type and Hurewicz-type results,~\rm Trans. A. M. S.~304 (1987), 431-467

\noindent [Ma]\ \ R. D. Mauldin,~\it One-to-one selections, marriage theorems,~\rm Amer. J. Math.~104 (1982), 823-828

\noindent [Mo]\ \ Y. N. Moschovakis,~\it Descriptive set theory,~\rm North-Holland, 1980

\noindent [P]\ \ T. C. Przymusinski,~\it On the notion of n-cardinality,~\rm Proc. Amer. Soc.~69 (1978), 333-338

\noindent [W]\ \ W. W. Wadge,~\it Thesis,~\rm Berkeley (1984) 

\end{document}